\documentclass[a4paper,11pt,reqno]{amsart}
\usepackage[latin1]{inputenc}
\usepackage[english]{babel}
\usepackage{amssymb,amsfonts, amsmath, color}
\usepackage[margin=2.7cm]{geometry}
\usepackage[table, dvipsnames]{xcolor}
\usepackage{graphicx, color, enumerate}
\usepackage[latin1]{inputenc}
\usepackage[active]{srcltx}
\usepackage{tikz}
\usepackage{pgf}
\usepackage{etex}
\usepackage{verbatim}
\usepackage{tikz-3dplot}
\usepackage{pgfkeys}
\usepackage{amsmath}
\usepackage{amsfonts}
\usepackage{amssymb}
\usepackage{amsthm}
\usepackage{algpseudocode}
\usepackage{float}
\usepackage{mathrsfs}
\usepackage{import}
\usepackage{geometry}
\usepackage{fancyhdr}
\usepackage{float}
\usepackage{fp}
\usepackage[colorlinks, citecolor=blue, linkcolor=red]{hyperref}

\newtheorem{theorem}{Theorem}[section]
\newtheorem{proposition}[theorem]{Proposition}
\newtheorem{corollary}[theorem]{Corollary}
\newtheorem{lemma}[theorem]{Lemma}
\newtheorem{remark}[theorem]{Remark}
\newtheorem{definition}[theorem]{Definition}
\newtheorem{example}[theorem]{Example}

\newcommand{\bcl}{\begin{center}}
\newcommand{\ecl}{\end{center}}
\newcommand{\brl}{\begin{right}}
\newcommand{\erl}{\end{right}}
\newcommand{\ben}{\begin{enumerate}}
\newcommand{\een}{\end{enumerate}}
\newcommand{\overliner}{\begin{array}}
\newcommand{\earr}{\end{array}}
\newcommand{\btab}{\begin{tabular}}
\newcommand{\etab}{\end{tabular}}
\newcommand{\bdoc}{\begin{document}}
\newcommand{\edoc}{\end{document}}
\newcommand{\beqy}{\begin{eqnarray}}
\newcommand{\eeqy}{\end{eqnarray}}

\newcommand{\beqi}{\begin{eqnarray*}}
\newcommand{\eeqi}{\end{eqnarray*}}
\newcommand{\bitem}{\begin{itemize}}

\newcommand{\eitem}{\end{itemize}}
\newcommand{\nln}{\newline}
\newcommand{\newt}{\newtheorem}

\newcommand{\pa}{\partial}
\newcommand{\re}{{I\!\!R}}
\newcommand{\Rn}{\R^N}
\newcommand{\xr}{x\in\R }
\newcommand{\x}{\times}
\newcommand{\dyle}{\displaystyle}
\newcommand{\ene}{{I\!\!N}}
\newcommand{\irn}{\int\limits_{\R^N}}
\newcommand{\io}{\int\limits_{\O}}
\newcommand{\meas}{{\rm meas\,}}
\newcommand{\dif}{\nabla_{xy}}
\newcommand{\sign}{{\rm sign}}
\newcommand{\map}{\longrightarrow }
\newcommand{\imp}{\Longrightarrow }
\renewcommand{\div}{\nabla\cdot }
\newcommand{\sen}{{\rm sen\,}}
\newcommand{\tg}{{\rm tg\,}}
\newcommand{\arcsen}{{\rm arcsen\,}}
\newcommand{\arctg}{{\rm arctg\,}}
\newcommand{\supp}{{\textsl supp\ }}
\newcommand{\ity}{\int_{-\iy}^{+\iy}}
\newcommand{\limit}{\lim\limits}
\newcommand{\limi}{\limit_{n\to\infty}}
\newcommand{\sumi}{\sum\limits_{n=1}^{\infty}}
\newcommand{\ulu}{\underline u}
\newcommand{\ulw}{\underline w}
\newcommand{\ulz}{\underline z}
\newcommand{\ulv}{\underline v}
\newcommand{\uls}{\underline s}
\newcommand{\olu}{\overline u}
\newcommand{\olv}{\overline v}
\newcommand{\ols}{\overline s}
\newcommand{\ob}{\overline\b}
\newcommand{\ovar}{\overline\var}
\newcommand{\wv}{\widetilde v}
\newcommand{\wu}{\widetilde u}
\newcommand{\ws}{\widetilde s}
\renewcommand{\a }{\alpha }
\renewcommand{\b }{\beta }
\newcommand{\g }{\gamma}
\newcommand{\G }{\Gamma }
\renewcommand{\d }{\delta }

\newcommand{\bd }{\bar d }
\newcommand{\D }{\Delta }
\newcommand{\e }{\varepsilon }
\newcommand{\z }{\zeta }
\renewcommand{\l }{\lambda }
\renewcommand{\L }{\Lambda }
\newcommand{\m }{\mu }
\newcommand{\n }{\nabla }
\newcommand{\s }{\sigma }
\newcommand{\Sig }{\Sigma }
\renewcommand{\t }{\tau }
\newcommand{\var }{\varphi }
\renewcommand{\o }{\omega }
\renewcommand{\O }{\Omega }
\newcommand{\R}{{\mathbb{R}}}
\newcommand{\bC}{{\bf C}}
\newcommand{\bZ}{{\bf Z}}
\newcommand{\bN}{{\bf N}}
\newcommand{\bQ}{{\bf Q}}
\newcommand{\bK}{{\bf K}}
\newcommand{\bI}{{\bf I}}
\newcommand{\bv}{{\bf v}}
\newcommand{\bV}{{\bf V}}
\newcommand{\LL}{\mathcal{L}}
\newcommand{\N}{\mathbb{N}}
\DeclareMathOperator{\suppo}{supp} \DeclareMathOperator{\di}{div}

\newenvironment{Proof}{\Rmovelastskip\vskip12pt
plus 1pt \noindent\em\rm}{\hfill {\qed \hskip .2cm}}

\begin{document}

\title[]{Phragm\`en-Lindel\"of type theorems \\ for parabolic equations \\ on infinite graphs}

\author{Stefano Biagi}

\address{\hbox{\parbox{5.7in}{\medskip \noindent{Stefano Biagi, \\Dipartimento di Matematica, \\Politecnico di Milano, \\Piazza Leonardo da Vinci 32, 20133, Milano, Italy \\ [3pt] \emph{E-mail address: }{\tt stefano.biagi@polimi.it}}}}}

\author{Giulia Meglioli}

\address{\hbox{\parbox{5.7in}{\medskip \noindent{Giulia Meglioli, \\Dipartimento di Matematica, \\Politecnico di Milano, \\Piazza Leonardo da Vinci 32, 20133, Milano, Italy \\ [3pt] \emph{E-mail address: }{\tt giulia.meglioli@polimi.it}}}}}

\author{Fabio Punzo}

\address{\hbox{\parbox{5.7in}{\medskip \noindent{Fabio Punzo, \\Dipartimento di Matematica, \\Politecnico di Milano, \\Piazza Leonardo da Vinci 32, 20133, Milano, Italy \\ [3pt] \emph{E-mail address: }{\tt fabio.punzo@polimi.it}}}}}



\begin{abstract} We obtain the Phragm\`en-Lindel\"of principle on combinatorial infinite weighted graphs for the Cauchy problem associated to a certain class of parabolic equations with a variable density. We show that the hypothesis made on the density is optimal.
\end{abstract}

\maketitle

\section{Introduction} \label{sec0}
We investigate uniqueness of possibly unbounded solutions to parabolic Cauchy problem of the following type:
\begin{equation}\label{problema}
 \begin{cases} \rho\,\partial_tu-\Delta u=f &\quad \text{ in } G\times(0,T]=:S_T \\
 u=u_0&\quad \text{ in } G\times\{0\}.
 \end{cases}
  \end{equation}
  Here, $(G, \omega, \mu)$ denotes an infinite graph equipped with edge weights $\omega$ and vertex measure $\mu$. The function $\rho > 0$ plays the role of a density, and $\Delta$ is the graph Laplacian. The initial data $u_0$ and the source term $f$ are prescribed.

The analysis of partial differential equations on graphs, particularly on infinite and weighted structures, has received significant attention in recent years (see, e.g., \cite{Grig2, KLW, Mu2}).
While elliptic equations have been widely explored (e.g., \cite{AS1, BP25, BiMePu, GLY1, GLY2, HW, MP2}), the parabolic setting has seen substantial development in works such as \cite{BCG, CGZ, EMM, EM, GT, HL, HMu, Huang, HKS, KR, LW2, LHN, Me, MP3, MPS2, Mu, SSV, Wu}.

This paper is devoted to establishing uniqueness results for solutions to \eqref{problema}, under appropriate growth conditions, even allowing for solutions that are not bounded. Our main approach relies on proving a Phragm\`en-Lindel\"of type principle for the problem in the graph setting (see Proposition \ref{prop1}, Theorems \ref{teo1}, \ref{teo3}). From this, uniqueness of solutions, possibly unbounded, follows as a direct consequence (see Corollaries \ref{cor2}, \ref{cor5}).

There exists a vast body of literature concerning uniqueness and Phragm\`en-Lindel\"of type results for parabolic equations in Euclidean space $\mathbb R^n$ (e.g., \cite{EKP, IKO, KKT, KP1, KP2, Meg, MP, MoPu1, MoPu2, PPT1, PPT2, SSV, Smi, Son}), as well as on Riemannian manifolds (e.g., \cite{BP1, BiMePuCCM, Grig, GrigHK, MR, Pinch, Pu1}). Our work extends this framework to the discrete and infinite setting of graphs. Some related results for elliptic equations on graphs are established in \cite{BP25} (see Remark \ref{paragone-ell}).

\subsection{Overview of our results}
We begin by formulating a general Phragm\`en-Lindel\"of principle (Proposition \ref{teo1}) under the assumption of an appropriate supersolution, which makes the result somewhat implicit. We then demonstrate that, for a large class of graphs, such supersolutions can be explicitly constructed when the density $\rho$ satisfies a decay condition that depends on a key geometric feature of the graph, known as the outer degree (or outer curvature). This leads to explicit uniqueness criteria (Theorems \ref{teo1}, \ref{teo2}).

On certain graph classes, particularly spherically symmetric trees, we verify that the decay assumptions on $\rho$ and the outer degree are optimal (Theorem \ref{prop71}, Corollaries \ref{cor1pa}). Indeed, when these conditions are violated, we can construct infinitely many bounded solutions, which directly implies non-uniqueness. The construction is nontrivial due to the absence of standard a priori estimates available in the Euclidean case, necessitating a tailored argument for the graph context.

Moreover, we show that on the integer lattice $\mathbb Z^n$, further uniqueness results can be obtained under faster decay of the density $\rho$ (see Theorem \ref{teo3}), and we prove that this threshold is sharp (Corollary \ref{cor3pa}). Similar results will be obtained also for anti-trees (see Corollary \ref{cora3Bis} and Theorem \ref{teo2antit}).

Finally, in Section \ref{Z2}, we show that in the special cases of  $\mathbb Z^2$, uniqueness follows without any constraint on the decay rate of $\rho$.
\medskip

We collect in the next table our main uniqueness results (see the forthcoming sections for the relevant notation).
\begin{table}[H]\label{table}
\centering
\begin{tabular}{||>{\columncolor[gray]{0.95}}c||c|c|c||}
\hline
\rowcolor[gray]{0.95}
\rule[-4mm]{0mm}{1cm}
 & Assumption on $\rho$ & Growth condition for $u$ & Optimality on $\rho$
\\
\hline
\rule[-8mm]{0mm}{2cm}
General $G$ & $\begin{gathered}\rho(x)\geq \frac{\mathfrak D_+(x)}{ r+1}\,e^{\rho_0\log^\beta(r+2)}\\[0.15cm] 
(0\leq\beta\leq 1)\end{gathered}$ & $e^{ B(r+1)\log^\beta(r+1)}$ & It depends on $G$
\\
\hline
\rule[-8mm]{0mm}{2cm}
$\mathbb Z^n$, $n\ge 3$ & $\begin{gathered} \rho(x)\geq \rho_0(1+|x|)^{-\alpha}\\[0.15cm] (0\leq\alpha\leq 2) \end{gathered}$ & $
\begin{cases}
e^{B|x|^{\min\{1, 2-\alpha\}}}, & \alpha\in [0, 2) \normalcolor\\
 e^{B\log^2(2+|x|^2)}, & \alpha =2
\end{cases}$ & Yes
 \\
\hline
\rule[-4mm]{0mm}{1cm}
$\mathbb Z^2$ & $\rho>0$ & $\log(\log(|x|^2+4))$ & Obvious
 \\
\hline
\rule[-4mm]{0mm}{1cm}
Tree  & $\rho(x)\geq\rho_0\, \frac{b}{ r+1} $ & $e^{ B(r+1)}$ & Yes
\\
\hline
\rule[-4mm]{0mm}{1cm}
Anti-tree & divergence of a certain series & it depends on $G$ and $\rho$ & Yes \\
\hline
\end{tabular}
\caption{An overview on our uniqueness results}
\end{table}

\normalsize

\smallskip
\subsection{Structure of the paper}
The paper is organized as follows. In Section \ref{mf} we provide the main definitions concerning the graph setting and the involved operators on graphs. Afterwards, in Section \ref{main} we state the main results: first the Phragm\'en-Lindel\"of principle and the uniqueness result, afterward non-uniqueness and optimality. Section \ref{Zn} is devoted to the case of the lattice $\mathbb Z^n$ which deserves a special attention since it differs from the general case. In Section \ref{auxiliary} we establish a weak maximum principle. The proof of the general Phragm\`en-Lindel\"of principle is given in Section \ref{proofPL}. Afterwards, in Section \ref{nonuniqueness} we construct proper solutions which demonstrate nonuniqueness and let us discuss optimality. Section \ref{integer} presents additional results specific to the lattice $\mathbb Z^n$. Finally, Sections \ref{antitree} and \ref{Z2} present further
developments on anti-trees and on the special case of \(\mathbb Z^2\),
respectively.

\section{Mathematical framework and the main result}\label{mf}\setcounter{equation}{0}

\subsection{The graph setting}
Let $G$ be a countably infinite set and $\mu:G\to (0,+\infty)$ be a measure on $G$ satisfying $\mu(\{x\})<+\infty$ for every $x\in G$ (so that $(G,\mu)$ becomes a measure space). Furthermore, let
\begin{equation*}
\omega:G\times G\to [0,+\infty)
\end{equation*}
be a symmetric, with zero diagonal and finite sum function, i.e.
\begin{equation*}
\begin{aligned}
&\text{(i)}\,\, \omega(x,y)=\omega(y,x)\quad &\text{for all}\,\,\, (x,y)\in G\times G;\\
&\text{(ii)}\,\, \omega(x,x)=0 \quad\quad\quad\,\, &\text{for all}\,\,\, x\in G;\\
&\text{(iii)}\,\, \displaystyle \sum_{y\in G} \omega(x,y)<\infty \quad &\text{for all}\,\,\, x\in G\,.
\end{aligned}
\end{equation*}
Thus, we define  \textit{weighted graph} the triplet $(G,\omega,\mu)$, where $\omega$ and $\mu$ are the so called \textit{edge weight} and \textit{node measure}, respectively. Observe that assumption $(ii)$ corresponds to ask that $G$ has no loops.
\smallskip

\noindent Let $x,y$ be two points in $G$; we say that
\begin{itemize}
\item $x$ is {\it connected} to $y$ and we write $x\sim y$, whenever $\omega(x,y)>0$;
\item the couple $(x,y)$ is an {\it edge} of the graph and the vertices $x,y$ are called the {\it endpoints} of the edge whenever $x\sim y$;
\item a collection of vertices $ \{x_k\}_{k=0}^n\subset G$ is a {\it path} if $x_k\sim x_{k+1}$ for all $k=0, \ldots, n-1.$
\end{itemize}

\noindent We are now ready to list some properties that the weighted graph $(G,\omega,\mu)$ may satisfy.
\begin{definition}\label{def01}
We say that the weighted graph $(G,\omega,\mu)$ is
\begin{itemize}
\item[(i)] {\em locally finite} if each vertex $x\in G$ has only finitely many $y\in G$ such that $x\sim y$;
\item[(ii)] {\em connected} if, for any two distinct vertices $x,y\in G$ there exists a path joining $x$ to $y$;
\end{itemize}
\end{definition}

\noindent For any $x\in G$, we define
\begin{itemize}
\item the {\it degree} of $x$ as $$\operatorname{deg}(x):=\sum_{y\in G}\omega(x,y);$$
\item the {\it weighted degree} of $x$ as 
\begin{equation} \label{eq:DefDegWMP}
 \operatorname{Deg}(x):=\frac{\operatorname{deg}(x)}{\mu(x)}.	
\end{equation}
\end{itemize}
Let now $d:G\times G\to[0,+\infty)$ be a \emph{distance on $G$}, that is,
\begin{itemize}
  \item[a)] $d(x,x)=0$ \quad \text{ for all } $x\in G$;
  \item[b)] $d(x,y)=d(y,x)$ \text{ for all } $x,y\in G$;
  \item[c)] $d(x,y)\leq d(x,z)+d(z,y)\quad \text{for all}\,\,\, x,y,z\in G$.
\end{itemize}
For any $x_0\in G$ and $r>0$ we define the ball $B_r(x_0)$ with respect to any metric $d$ as
\[B_r(x_0):=\{x\in G\,:\, d(x,x_0)<r\}\,.\]
Furthermore, we define the \textit{jump size} $s>0$ of a pseudo metric $d$ as
\begin{equation}\label{e14f}
s:=\sup\{d(x,y) \,:\, x,y\in G, \omega(x,y)>0\}.
\end{equation}
For a more detailed understanding of the objects introduced so far, we refer the reader to \cite{GT,HKMW,HKS,MP2}.
 In this paper, we always make the following assumptions
\begin{equation}\label{e7f}
\begin{split}
(G, \omega, \mu) \text{ is a connected, locally finite, weighted graph}.
\end{split}
\end{equation}

\subsection{Difference and Laplace operators} Let $\mathfrak F$ denote the set of all functions $f: G\to \mathbb R$. 	
	For any $f\in \mathfrak F$ and for all $x,y\in G$, let us give the following
\begin{definition}\label{def1}
Let $(G, \omega,\mu)$ be a weighted graph. For any $f\in \mathfrak F$,
\begin{itemize}
\item the {\em difference operator} is
\begin{equation*}
\nabla_{xy} f:= f(y)-f(x)\,;
\end{equation*}
\item the {\em (weighted) Laplace operator} on $(G, \omega, \mu)$ is
\begin{equation}\label{eq11}
\Delta f(x):=\frac{1}{\mu(x)}\sum_{y\in G}[f(y)-f(x)]\omega(x,y)\quad \text{ for all }\, x\in G\,.
\end{equation}
\end{itemize}
\end{definition}
Clearly,
\[\Delta f(x)=\frac 1{\mu(x)}\sum_{y\in G}(\dif f)\omega(x,y)\quad \text{ for all } x\in G\,.\]
It is straightforward to show, for any $f,g\in \mathfrak F$, the validity of
\begin{itemize}
\item  the {\it product rule}
\begin{equation*}
\nabla_{xy}(fg)=f(x) (\nabla_{xy} g) + (\nabla_{xy} f)g(y) \quad \text{ for all } x,y\in G\,;
\end{equation*}
\item the {\it integration by parts formula}
\begin{equation}\label{e4f}
\sum_{x\in G}[\Delta f(x)] g(x) \mu(x)=-\frac 1 2\sum_{x,y\in G}(\dif f)(\dif g) \omega(x,y)\,,
\end{equation}
provided that at least one of the functions $f, g\in \mathfrak F$ has {\it finite} support.
\end{itemize}



\subsection{Outer and inner degrees}\label{degree}
We introduce some basic definitions following \cite[Chapter 9]{KLW}.

\noindent We racall that the {\it combinatorial graph distance} on $G$, is the distance which, for any two vertices $x, y\in G$, counts the least number of edges in a path between $x$ and $y$; we name it $\bd$.

\noindent Let $\Omega\subset G$ be finite subset. Define the distance from any $x\in G$ to the subset $\Omega$
\[\bd(x, \Omega):=\min_{y\in \Omega} \bd(x, y)\quad \forall x\in G\,.\]
With an abuse of notation we write $\bd(x, y)$ to indicate the distance between any two points $x, y\in G$, and $\bd(x, \Omega)$ to denote the distance from the point $x\in G$ to the set $\Omega\subset G$.

\noindent For any $r\in \mathbb N_0$, let
\[\mathcal S_r(\Omega):=\{x\in G\,: \bd(x, \Omega)=r\,\}.\]
Given $f\in \mathfrak F$, we say that $f$ is {\it spherically symmetric w.r.t.} $\Omega$ if
\[f(x)=f(y) \quad \text{ whenever } \bd(x, \Omega)=\bd(y, \Omega).\]
In this case, with a slight abuse of notation, we write
\[f(x)=f(r) \quad \forall\,\,x\in \mathcal S_r(\Omega)\,.\]
For any $x\in G$ with $r\equiv r(x):=\bd(x, \Omega)\geq 1$, let
\[\mathfrak D_+(x):=\frac1{\mu(x)}\sum_{y\in \mathcal S_{r+1}(\Omega)}\omega(x,y), \quad \mathfrak D_-(x):=\frac1{\mu(x)}\sum_{y\in \mathcal S_{r-1}(\Omega)}\omega(x,y)\,.\]
The function $\mathfrak D_+:G\to [0, +\infty)$ is called {\it outer degree} (or {\it outer curvature}) w.r.t. $\Omega$, whereas $\mathfrak D_-:G\to [0, +\infty)$ is called {\it inner degree} (or {\it inner curvature}) w.r.t. $\Omega$, (see \cite{AS1}).

The weighted graph $(G, \mu, \omega)$, endowed with the combinatorial distance $r$, is said to be
{\em weakly spherically symmetric} with
respect to a finite subset $\Omega\subset G$, if the outer and inner degrees $\mathfrak D_{\pm}$ are spherically symmetric with respect to $\Omega$. 

Therefore, on a weakly symmetric graph,
\[\mathfrak D_{\pm}(x)=\mathfrak D_{\pm}(r)\quad 
\forall\,\,x\in S_r(\Omega)\,.\]
\begin{remark} \label{rem:SferaOmegafinita}
	Let assumption \eqref{e7f} be in force, and let
	$\bar{d}$ the combinatorial distance on $G$
	introduced above. If $\Omega\subseteq G$
	is a \emph{finite set}, it follows from the 
	\emph{local finiteness} of $G$ that
	$$\Omega_R = \{x\in G:\,\bar{d}(x,\Omega) < R\}$$
	is a finite set for every $R > 0$. We will repeatedly exploit
	this fact in the sequel.
\end{remark}

\section{Main results} \label{main}\setcounter{equation}{0}

\noindent We have already stated in \eqref{e7f}
the main hypotheses on the weighted graph $(G,\omega,\mu)$.
In order to state our main results, we now fix
the precise definition
of sub/supersolution to problem
\eqref{problema} (and, more generally, to
the Cauchy-Dirichlet problem associated
with $\rho\,\partial_t-\Delta$).
\vspace{0.1cm}

Throughout the rest of the paper,
we exploit the following notation.
\begin{enumerate}
	\item[(N1)] We set 
$
\LL:=\rho\,\partial_t-\Delta.
$
\vspace{0.1cm}

\item[(N2)] If $\Omega\subseteq G$ is an arbitrary set, we indicate by $\mathfrak{F}(\Omega)$ the set of all functions $f:\Omega\to\R$. In the particular case $\Omega=G$, we simply write $\mathfrak F$,
 consistently with the notation introduced in the previous section.
\vspace{0.1cm}

\item[(N3)] Let $T\in(0,+\infty)$, and let $\Omega\subseteq G$ be an arbitrary set. We set
$$\Omega^T = \Omega\times (0,T];
\qquad
\Omega^T_0 = \Omega\times[0,T].$$
Accordingly, we indicate by $\mathfrak{F}^T(\Omega)$ 
(resp., $\mathfrak{F}^T_0(\Omega)$)
the set of all functions $f:\Omega^T\to\R$ (resp., $f:\Omega^T_0\to\R$). In the particular case when
$\Omega = G,$
we simply write $$
\text{$\mathfrak{F}^T$\,\,(resp., $\mathfrak{F}^T_0$)}.$$
Finally, we define \(\mathfrak R^T_0(\Omega)\) as the
set of all functions \(f\in\mathfrak F^T_0\), i.e., 
\emph{defined on the whole of}
\(G\times\overline [0,T]\), satisfying the following 
\emph{regularity properties}:
$$u(x,\cdot)\in C([0,T])\cap C^1((0,T])\quad \forall\,\,x\in\Omega.$$
Hence, the dependence on $\Omega$ of the space $\mathfrak R^T_0(\Omega)$
lies in the fact that we prescribe the time regularity of $f(x,\cdot)$
only for $x\in\Omega$.
If $\Omega=G$, we simply write $\mathfrak R^T_0$.
	\end{enumerate}
\begin{definition}\label{defsol}
Given $T\in(0,+\infty)$, $f\in\mathfrak{F}^{T}$ and $u_0\in\mathfrak F$, we say
 that a function $u\in\mathfrak{F}_0^T$ is a subsolution \emph{[}resp.\,su\-per\-so\-lu\-tion\emph{]} to problem \eqref{problema} if
 \begin{itemize}
  \item[i)] $u\in\mathfrak{R}_0^T$ (see notation \emph{(N3)});
 \item[ii)] $u$ solves the inequality $\LL u \le [\ge] \,\,f$  in $G^T$;
   \item[iii)] $u(x,0)\le [\ge] \,\,u_0(x)$ pointwise in $G$.
\end{itemize}
Moreover, we say that $u$ is a solution of  \eqref{problema} if it is both a subsolution and a supersolution.
\end{definition}

\begin{definition}\label{defsolOmega}
Let $\Omega\subsetneq G$ be 
a \emph{proper subset}, and let
$T\in(0,+\infty)$. Given any $f\in \mathfrak{F}^T(\Omega)$, 
any $g\in \mathfrak{F}^T_0(G\setminus\Omega)$ and 
any $u_0\in\mathfrak{F}(\Omega)$, we say that a function 
$$u\in\mathfrak{F}^T_0$$
is a subsolution \emph{[}resp.\,supersolution\emph{]} of the $\LL$\,-\,Dirichlet problem
 \begin{equation} \label{eq:pbDir}
  \begin{cases}
 \LL u = f & \text{in $\Omega^T$} \\
 u = g & \text{in $(G\setminus\Omega)^T_0$}\\
 u=u_0& \text{in $\Omega\times\{0\}$},
 \end{cases}
 \end{equation}
 if the following conditions hold:
 \begin{itemize}
   \item[i)] $u\in\mathfrak{R}^T_0(\Omega)$;
 \item[ii)] $u$ solves the inequality $\LL u \le [\ge] \,\,f$  in $\Omega^T$;
   \item[iii)] $u(x,t) \le [\ge] \,\,g(x,t)$ pointwise in $(G\setminus\Omega)^T_0$;
   \item[iv)] $u(x,0)\le [\ge] \,\,u_0(x)$ pointwise in $\Omega$.
  \end{itemize}
Finally, we say that $u$ is a solution of problem \eqref{eq:pbDir} if $u$ is both  a subsolution and a supersolution of this problem.
 \end{definition}

\subsection{Phragm\`en-Lindel\"of principle and uniqueness results}
The first main result
of this paper is a general Phrag\-m\`en\--Lin\-del\"of-type principle, which reads as follows. 

Throughout the rest of the paper, we tacitly understand that
\begin{center}
\emph{$T$ is a fixed number in $(0,+\infty)$}.
\end{center}
\begin{proposition}\label{prop1}
Let assumption \eqref{e7f} be satisfied, and let $x_0\in G$ . Moreover, let $\rho\in \mathfrak F$
be such that $\rho>0$ on $G$. Suppose that there exists $Z\in \mathfrak{R}^T_0$, $Z>0$ in 
	$G^T_0$, such that
\begin{equation}\label{e1p}
\rho(x)\,\partial_t Z(x,t)-\Delta Z(x,t)\ge0\quad\text{for all}\,\,\,(x,t)\in G^T\,.
\end{equation}
Let $u$ be a subsolution of problem \eqref{problema} with $f\equiv0$, $u_0\equiv0$ fulfilling
\begin{equation}\label{e3p}
\limsup_{d(x,x_0)\to +\infty}\left\{\max_{t\in[0,T]} \frac{u(x,t)}{Z(x,t)}\right\}\le 0\,.
\end{equation}
Then
\[u\le 0 \quad \text{ in } G^T_0\,.\]
\end{proposition}
\begin{remark} \label{rem:Equivalenteuplus}
	Let the assumptions and the notation of Proposition \ref{prop1} be in force. We explicitly stress that condition \eqref{e3p} is actually equivalent to the following one
	$$\lim_{d(x,x_0)\to+\infty}\left\{\max_{t\in[0,T]} \frac{u^+(x,t)}{Z(x,t)}\right\}= 0,$$
	where $u^+ = \max\{u,0\}$ is the positive part of $u$.
\end{remark}
The previous proposition is an abstract result: it shows that a
Phragm\`en--Lindel\"of principle fol\-lows once a suitable positive supersolution
\(Z\) of \eqref{e1p} is available. The next results are devoted to showing that,
under explicit assumptions on the density \(\rho\) and on the geometry of the
graph, such supersolutions can indeed be constructed.

To this end, throughout the sequel, whenever \(\Omega\subseteq G\) is finite, we
write
\begin{equation}\label{def:r}
r\equiv r(x):= \bar d(x, \Omega)\quad \forall\, x\in G\,.
\end{equation}
where \(\bar d\) is the combinatorial distance.

\begin{theorem}\label{teo1}
Let assumption \eqref{e7f} be satisfied, and let $\Omega\subset G$ be a finite set.  Moreover, let
$\rho\in\mathfrak{F}$ be such that
\begin{equation}\label{e12f}
\rho(x)\geq\rho_0 \frac{\mathfrak D_+(x)}{ r+1} \quad \textrm{for all}\;\; x\in G,
\end{equation}
for some
$\rho_0>0.$
 Let $u$ be a subsolution of problem \eqref{problema} with $f\equiv u_0\equiv0$ fulfilling
\begin{equation}\label{e16p}
\limsup_{r\to +\infty}\frac1{\tilde Z(x)}\left\{\max_{t\in[0,T]}{u(x,t)}\right\} \le 0\,,
\end{equation}
where for some $B>0$,
\begin{equation}\label{e40f}
  \tilde Z(x):= e^{ B(r+1)}\,,\quad\text{for all}\,\,\,x\in G.
\end{equation}
Then
\[u \leq 0 \quad\text{in}\,\,\, G^T_0.\]
\end{theorem}

\begin{theorem}\label{teo2}
Let assumption \eqref{e7f} be satisfied. Let $\Omega\subset G$ be a finite subset. Suppose that
\begin{equation}\label{eq:rho-gen}
\rho\in \mathfrak F, \quad \rho(x)\geq \frac{\mathfrak D_+(x)}{r+1}\,e^{\rho_0\log^\beta(r+2)} \quad \textrm{for all}\;\; x\in G,
\end{equation}
for some $\beta\in (0, 1]$ and $\rho_0>0$.
 Let $u$ be a subsolution of problem \eqref{problema} with $f\equiv u_0\equiv0$ fulfilling
\begin{equation}\label{eq:limite}
\limsup_{r\to +\infty}\frac1{\hat Z(x)}\left\{\max_{t\in[0,T]}{u(x,t)}\right\} \le 0\,,
\end{equation}
where, for some \(B\in\left(0,\frac{\rho_0}2\right)\),
\begin{equation}\label{eq:soprasol}
  \hat Z(x):= e^{ B(r+1)\log^\beta(r+2)}\,,\quad\text{for all}\,\,\,x\in G.
\end{equation}
Then
\[u \leq 0 \quad\text{in}\,\,\, G^T_0.\]
\end{theorem}

We can immediately infer the following uniqueness results.

\begin{corollary}\label{cor2}
Let assumption \eqref{e7f} be satisfied. Let $\rho\in \mathfrak F$, $\rho> 0$ and let $x_0\in G$ be some reference point. Suppose that there exists $Z\in \mathfrak {R}^T_0$, $Z>0$ in $ G^T_0$ such that \eqref{e1p} holds. Then there exists at most one solution to problem \eqref{problema} such that
\begin{equation}\label{el3Bis}
\lim_{d(x,x_0)\to +\infty}\left\{\max_{t\in[0,T]} \frac{|u(x,t)|}{Z(x,t)}\right\}=0\,.
\end{equation}
\end{corollary}

\begin{corollary}\label{cor3}
Let assumption \eqref{e7f} be satisfied and assume \eqref{e12f}. Let $\tilde Z$ be as defined in \eqref{e40f}. Then there exists at most one solution to problem \eqref{problema} such that
\[\lim_{r\to +\infty}\frac1{\tilde Z(x)}\left\{\max_{t\in[0,T]}{|u(x,t)|}\right\} =0\,.\]
\end{corollary}

\begin{corollary}\label{cor4}
Let assumption \eqref{e7f} be satisfied and assume \eqref{eq:rho-gen}. Let $\hat Z$ be as defined in \eqref{eq:soprasol}. Then there exists at most one solution to problem \eqref{problema} such that
\[\lim_{r\to +\infty}\frac1{\hat Z(x)}\left\{\max_{t\in[0,T]}{|u(x,t)|}\right\} =0\,.\]
\end{corollary}

\begin{remark}\label{paragone-ell}
Let $Z\in \mathfrak F$ be such that
$$
\Delta Z(x)\le \rho(x)\quad \text{for all}\,\,\,x\in G,
$$
and $\displaystyle \inf_{x\in G} Z(x)>0$.
In view of Lemma \ref{lemmaelliptic}, Proposition \ref{prop1} can be applied with
\begin{equation*}
\limsup_{d(x, x_0)\to +\infty} \frac{1}{Z(x)}\left\{\max_{t\in[0,T]} u(x,t)\right\}\leq 0\,
\end{equation*}
instead of \eqref{e3p}. 
In addition, Corollary \ref{cor2} holds with \eqref{el3Bis} replaced by 
\[\lim_{d(x,x_0)\to +\infty}\left\{\max_{t\in[0,T]} \frac{|u(x,t)|}{Z(x)}\right\}=0\,.
\]
In \cite{BP25}, certain supersolutions $Z$ of \eqref{el1} are constructed. As noted above, such supersolutions are expected to yield results analogous to Theorems \ref{teo1}, \ref{teo2}, and \ref{teo3}, as well as Corollaries \ref{cor3}, \ref{cor4}, and \ref{cor5}, albeit under different growth conditions at infinity. In contrast, in the present paper we construct supersolutions that explicitly depend on time. This allows us to establish a Phragm\`en-Lindel\"of principle under significantly weaker growth restrictions at infinity for the solution $u$. As a consequence, much larger uniqueness classes of solutions are obtained.
\end{remark}

\smallskip

\begin{remark}\label{ossmeas}
Observe that problem \eqref{problema} with $f=0$ and $u_0=0$ can be equivalently rewritten as
\begin{equation}\label{enf15}
    \begin{cases}
        u_t=\tilde \Delta u 
        & \text{in } \tilde G\times(0,+\infty),\\
        u(x,0)=0
        & \text{in } \tilde G,
    \end{cases}
\end{equation}
where \(\tilde G=(G,\omega,\tilde \mu)\), 
\[
    \tilde \mu(x)=\rho(x)\mu(x),
\]
and \(\tilde \Delta\) denotes the corresponding weighted Laplacian.

It is worth observing, however, that this reformulation does not lead
to any substantial simplification for the purposes of the present
paper.

Uniqueness for problem \eqref{enf15} has been extensively investigated in
the literature. Let us stress that our approach yields Phragm\'en--Lindel\"of
type principles, rather than uniqueness results alone. The known uniqueness
results mainly concern \emph{bounded} solutions and are closely related to
stochastic completeness; see, for instance, \cite{KLW}.
In other works (see e.g. \cite{Huang, HKS, KLW, Me}), uniqueness is proved for solutions satisfying suitable
integral growth conditions, which are tied to the measure considered on the
graph and to an intrinsic distance. 

By contrast, in the present paper we impose direct pointwise conditions
at infinity, which are independent of the measure; moreover, we do not require the distance to be intrinsic. More precisely, our
assumptions are formulated in terms of the behaviour of \(\rho(x)\) and of
the pointwise growth of the solution. In particular, our approach allows us
to treat solutions that may be unbounded.

Let us also note that pointwise growth conditions at infinity and integral
growth assumptions are not comparable in general, unless further information
on the measure of the graph is available. This is because integral conditions
also reflect the volume growth associated with the underlying measure, while
pointwise conditions do not.
\end{remark}

Now, we consider a special kind of weakly symmetric graphs, the so called \textit{spherically symmetric trees}. Let $(G, \omega, \mu)$ be a weakly symmetric graph w.r.t. $\Omega=\{o\}$, for some
fixed point $o\in G$
(which is usually referred to as the \emph{root of $G$}).
Suppose that
\begin{itemize}
\item $\omega:G\times G\to \{0, 1\}$;
\vspace*{0.1cm}

\item $\omega|_{S_m(\Omega)\times S_m(\Omega)}=0$;
\vspace*{0.1cm}

\item $\mu(x)=1 \text{ for every } x\in G;$
\vspace*{0.1cm}

\item there exists $b:\mathbb N\to\mathbb{N}$, which is called the {\it branching function}, such that
\[\mathfrak D_+(x)=b(m), \quad \mathfrak D_-(x)=1\quad \text{for every $x\in S_m(\Omega)$ and
$m\in \mathbb N$}.\]
\end{itemize}

From Corollary \ref{cor3} we immediately deduce the following result. 

\begin{corollary}\label{cortree}
Let  $(G, \omega, \mu)$ be a spherically symmetric tree as above, with
  \emph{constant branching function} $b(r) = b_0\geq 2$. Suppose that
\begin{equation}\label{enf16}
\rho(x)\geq\rho_0 \frac{\mathfrak D_+(x)}{ r+1} \quad \textrm{for all}\;\; x\in G,
\end{equation}
for some $\rho_0>0$. Let $\tilde Z$ be as defined in \eqref{e40f}. Then there exists at most one solution to problem \eqref{problema} such that
\[\lim_{r\to +\infty}\frac1{\tilde Z(x)}\left\{\max_{t\in[0,T]}{|u(x,t)|}\right\} =0\,.\]
\end{corollary}

On a spherically symmetric tree, in the special case that $\rho$ is radial we can discuss the previous result also by means of Remark \ref{ossmeas}.

\begin{example}[The homogeneous tree and the modified measure]\label{treesc}
Let  $(G, \omega, \mu)$ be a spherically symmetric tree as above, with
  \emph{constant branching function} $b(r) = b_0\geq 2$.
  
  Assume that \(\rho\) satisfies \eqref{enf16}
  and \emph{is radial}, namely
\[
\rho(x)=\rho_r
\qquad\text{if } d(x, o)=r.
\]
The equation
\[
\rho(x)u_t=\Delta u
\]
can be rewritten as a standard heat equation
\[
u_t=\widetilde\Delta u
\]
on the same weighted graph, but with respect to the modified vertex measure
\[
\widetilde\mu(x):=\rho(x)\mu(x).
\]
Notice that, since \(\mu\equiv1\), we have
\[
\widetilde\mu(x)=\rho_r
\qquad\text{if } d(x,o)=r.
\]
The transformed graph is again spherically symmetric. Its weighted volume
balls are given by
\[
\widetilde V(R)
:=
\sum_{x\in B_R}\widetilde\mu(x)
=
\sum_{k=0}^R |S_k|\rho_k,
\]
where
\[
S_k:=\{x\in G:\ |x|=k\}.
\]
Moreover, for \(k\ge1\),
\[
|S_k|=(b_0+1)b_0^{k-1},
\]
and the weighted boundary of \(B_R\) is
\[
\partial B_R
:=
\sum_{\substack{x\in S_R\\ y\in S_{R+1}}}\omega(x,y)
=
b_0 |S_R|.
\]
By the standard stochastic completeness criterion for spherically symmetric
trees (see e.g. \cite[Chapter 9]{KLW}), the transformed graph \((G,\omega,\widetilde\mu)\) is stochastically
complete if and only if
\[
\sum_{R=1}^{+\infty}
\frac{\widetilde V(R)}{\partial B_R}
=
+\infty.
\]
We now verify this condition. Since
\[
\widetilde V(R)
=
\sum_{k=0}^R |S_k|\rho_k
\ge
|S_R|\rho_R,
\]
we get, using our assumption \eqref{enf16} on \(\rho_R\),
\[
\frac{\widetilde V(R)}{\partial B_R}
\ge
\frac{|S_R|\rho_R}{b_0|S_R|}
=
\frac{\rho_R}{b_0}
\ge
\frac{\rho_0}{R+1}.
\]
Therefore
\[
\sum_{R=1}^{+\infty}
\frac{\widetilde V(R)}{\partial B_R}
\ge
\rho_0
\sum_{R=1}^{+\infty}
\frac1{R+1}
=
+\infty.
\]
Hence the graph with modified measure
\[
(G,\omega,\widetilde\mu),
\qquad
\widetilde\mu=\rho\mu,
\]
is stochastically complete. Therefore (\cite{KLW}) we have uniqueness of {\em bounded solutions} to problem \eqref{enf15}.

\smallskip

This shows that, on the homogeneous tree, our pointwise assumption on
\(\rho\) implies stochastic completeness of the graph obtained by absorbing
\(\rho\) into the measure. In this sense, the result is consistent with the
classical measure-theoretic approach. However, our theorem gives more than
bounded uniqueness: it yields a Phragm\'en--Lindel\"of type principle under
explicit pointwise growth conditions at infinity. Moreover, in general we do not require that $\rho$ is radial. 
\end{example}

\begin{example}[The homogeneous tree: comparison of uniqueness classes]
(a) Let \(G\) be the homogeneous tree with branching function
\(b(r)=b_0\ge2\). Let \(\rho\equiv1\). The structural assumption
\[
\rho(x)\ge \rho_0\frac{\mathfrak{D}_+(x)}{r+1}
\]
is satisfied, for instance, by choosing \(0<\rho_0\le \frac1{b_0}\), since
\[
1=\rho(x)\ge \rho_0\frac{b_0}{r+1}
\qquad\text{for all } r\ge0.
\]
Our Phragm\'en--Lindel\"of principle in Theorem \ref{teo1} then gives uniqueness in the pointwise
class
\[
\limsup_{r\to+\infty}
\frac{\max_{t\in[0,T]}u(x,t)}
{\exp(B(r+1))}
\le0
\]
for some \(B>0\). In particular, it allows solutions with exponential
pointwise growth in the combinatorial distance.

Let us compare this with the integral uniqueness class of
\cite[Theorem 1.3]{HKS}, \cite[Theorem 3.1]{Me}. Since the graph has bounded degree and
\(\mu\equiv1\), the combinatorial distance is equivalent to an intrinsic
metric. Moreover, the graph is globally local for functions of the form
\[
f(R)=aR,
\qquad a>0.
\]
The corresponding integral growth condition has the form
\[
\int_0^T
\sum_{x\in B(o,R)}
|u(x,t)|^2\,d\mu(x)\,dt
\le
C e^{aR},
\]
for every \(R>0\), up to harmless changes of the constants.

On the homogeneous tree,
\[
|S_0|=1,
\qquad
|S_r|=(b_0+1)b_0^{r-1}
\quad\text{for } r\ge1.
\]
Hence
\[
|S_r|\asymp b_0^r.
\]
If \(u\) is radially controlled, namely
\[
|u(x,t)|\le U_r
\qquad\text{whenever } |x|=r,
\]
then the contribution of the sphere \(S_r\) to the above integral condition is
of order
\[
|S_r|U_r^2
\asymp
b_0^r U_r^2.
\]
Thus the condition
\[
\int_0^T
\sum_{x\in B(o,R)}
|u(x,t)|^2\,d\mu(x)\,dt
\le
C e^{aR}
\]
is compatible, at the level \(r=R\), with
\[
b_0^r U_r^2\lesssim e^{ar}.
\]
Equivalently,
\[
U_r
\lesssim
\exp\left(\frac{a-\log b_0}{2}\,r\right).
\]
Since \(a>0\) is arbitrary, this gives an exponential uniqueness class in the
combinatorial distance. In particular, for every fixed \(B>0\), choosing
\[
a>2B+\log b_0
\]
shows that the pointwise bound
\[
|u(x,t)|\le C e^{B(r+1)}
\]
is compatible with the integral growth condition in \cite{HKS}.

Thus, in the homogeneous tree with \(\rho\equiv1\), the two approaches lead to
comparable uniqueness classes: our condition is pointwise and explicit, while
the condition in \cite{HKS} is integral and includes the exponential
volume growth of the spheres.

(b) In principle, more general weights \(\rho(x)\) could also be considered. 
However, the application of the results in \cite{HKS} would then no longer be immediate, 
because one would have to take into account both the intrinsic distance associated 
with the modified measure and the volume growth of balls with respect to the new measure. 
In general, this is a nontrivial matter.

\end{example}

\subsection{Optimality and nonuniqueness results} The main aim of this section is to provide
a general sufficient
condition for the existence of \emph{infinitely many solutions}
of problem \eqref{problema}; as we will see, thanks to this result we are
able to show that our uniqueness in Theorem \ref{teo1} is \emph{optimal}.
\vspace*{0.1cm}

To state the results of this section, we need to require some
additional assumptions on the graph $G$; more precisely,
together with assumption \eqref{e7f} we assume that
\begin{equation} \label{eq:extraassump}
 \begin{split}
  (i)& \,\,\text{there exists a \textit{pseudo metric}}\,\, d \,\,\,
\text{such that the jump size $s$ is finite}; \\
(ii)&\,\,\text{the ball}\,\,\, B_r(x) \,\,\,\text{with respect to}\,\,\, d\,\,\, \text{is a finite set, for any}\,\,\, x\in G,\,\,\, r>0.
 \end{split}
\end{equation}
\begin{remark} \label{rem:CasoParticoleHySopra}
 We explicitly stress that,  in view of Remark \ref{rem:SferaOmegafinita}, the above assumption 
 \eqref{eq:extraassump} is \emph{always satisfied}
 by any graph $G$ fulfilling \eqref{e1p}, provided that
 $$d = \bar{d}$$
 is the combinatorial distance defined above.
\end{remark}

\begin{theorem}\label{prop71}
Let assumptions \eqref{e7f}\,-\,\eqref{eq:extraassump} be in force and let $\rho\in\mathfrak{F}$
be such that $\rho>0$ on $G$. Moreover, let $u_0\in\mathfrak{F}$ be a \emph{non-negative and bounded} function on $G$.

We assume that there exist a function $h\in\mathfrak{F}$ and a ball $B_{\hat{R}}(o)\subseteq G$ such that
\begin{equation} \label{eq710}
\begin{split}
\mathrm{i)}&\,\,\Delta h\leq -\rho \quad \text{in}\,\,\, G\setminus B_{\hat R}(o),\\
\mathrm{ii)}&\,\,h>0\quad \text{in}\,\,\,G, \\
\mathrm{iii)}&\,\,\text{$h(x)\to 0$ as $d(x,o)\to+\infty$}.
\end{split}
\end{equation}
Then there exist infinitely many bounded solutions $u$ of problem \eqref{problema}. In
particular, for every fixed $\gamma\in \R$ with
\begin{equation} \label{eq:assumptionOnuzero} 
\gamma\geq M_0= \sup_G u_0\in[0,+\infty),	
\end{equation}
there exists a solution $u$ to problem \eqref{problema} with $f\equiv0$ such that
\begin{equation} \label{eq:limituNonUnique}
\text{$u(x,t_0)\to\gamma$ as $d(x,o)\to+\infty$}\quad\text{for every $t_0 \in (0,T)$}.
\end{equation}
\end{theorem}

Now we show that, on trees, the results in Theorem \ref{teo1} and
Corollary \ref{cor3} are sharp. More precisely, we show that if condition \eqref{e12f} fail, then Theorem \ref{prop71} can be applied, therefore
infinitely many bounded solutions of problem \eqref{problema} exist.
\vspace*{0.1cm}

From Theorem \ref{prop71}, after having exhibited the requested barrier $h$, we will deduce the following consequences.

\begin{corollary}\label{cor1pa}
Let  $(G, \omega, \mu)$ be a spherically symmetric tree as above, with
  \emph{constant branching function} $b(r) = b_0\geq 2$.
  Moreover, let
  $u_0\in\mathfrak{F}$ be a \emph{non-negative and bounded} function on $G$.
    Assume that $\rho\in\mathfrak{F},\,\rho> 0$
on $G$ fulfills
\[\rho(x) \leq c_0 (1+ r)^{-\alpha}\quad \text{ for any } x\in G,\]
for some $c_0>0, \alpha>1$. Then for every fixed $\gamma\in \R$ with satisfying
\eqref{eq:assumptionOnuzero}
there exists a solution $u$ to problem \eqref{problema} with $f\equiv0$ satisfying \eqref{eq:limituNonUnique}.
\end{corollary}

\begin{example}[Homogeneous tree and stochastic incompleteness of the modified graph]\label{treesinc}
Let \(G\) be the homogeneous tree with branching function
\(b(r)=b_0\ge2\). Let \(\rho\) be radial and assume, for instance, that
\[
\rho(x)=\rho_r:=(1+r)^{-\alpha}
\qquad\text{if } d(x, o)=r,
\]
with
\[
\alpha>1.
\]
Then \(\rho\) satisfies
\[
\rho(x)\le c_0(1+r)^{-\alpha}
\]
for a suitable \(c_0>0\), as in Corollary \ref{cor1pa}.
Let us rewrite
\[
\rho(x)u_t=\Delta u
\]
as the standard heat equation
\[
u_t=\widetilde\Delta u
\]
on the same weighted tree, but with respect to the modified measure
\[
\widetilde\mu(x):=\rho(x)\mu(x).
\]
Since \(\mu\equiv1\), we have
\[
\widetilde\mu(x)=\rho_r=(1+r)^{-\alpha}
\qquad\text{if } d(x, o)=r.
\]

We now verify that the weighted graph
\[
(G,\omega,\widetilde\mu)
\]
is stochastically incomplete. Since the tree is homogeneous and \(\rho\) is
radial, the modified graph is still spherically symmetric. Therefore we may
use the stochastic completeness criterion for weakly spherically symmetric
graphs (see \cite[Chapter 9]{KLW}). Stochastic
incompleteness is equivalent to
\[
\sum_{R=1}^{+\infty}
\frac{\widetilde V(R)}{\partial B_R}
<
+\infty,
\]
with the same notation as in Example \ref{treesc}.
Observe that
\[
\widetilde V(R)
=
\sum_{k=0}^R |S_k|\rho_k
=
1+\sum_{k=1}^R (b_0+1)b_0^{k-1}(1+k)^{-\alpha}.
\]
Since the last term dominates the geometric sum, there exists a constant
\(C>0\) such that
\[
\widetilde V(R)
\le
C\, b_0^R(1+R)^{-\alpha}.
\]
Consequently,
\[
\frac{\widetilde V(R)}{\partial B_R}
\le
C\,\frac{b_0^R(1+R)^{-\alpha}}{(b_0+1)b_0^R}
\le
C(1+R)^{-\alpha}.
\]
Since \(\alpha>1\), we get
\[
\sum_{R=1}^{+\infty}
\frac{\widetilde V(R)}{\partial B_R}
<+\infty.
\]
Hence, by the above criterion, the modified graph
\[
(G,\omega,\widetilde\mu)
\]
is stochastically incomplete. Therefore the solution to problem \eqref{enf15} is not unique, in the class of bounded functions.
Thus, in this example, the decay condition
\[
\rho(x)\le c_0(1+r)^{-\alpha},
\qquad \alpha>1,
\]
implies, after absorbing \(\rho\) into the vertex measure, the stochastic
incompleteness of the transformed graph. This is consistent with
Corollary \ref{cor1pa}, which yields nonuniqueness for the Cauchy
problem.

However, we have to observe that in Corollary \ref{cor1pa} we do not require that $\rho$ is radial; furthermore, we prove the existence of infinitely many bounded solutions to problem \eqref{enf15} starting from the same initial given condition $u_0\equiv 0$. This is a stronger conclusion than the previous one. 
\end{example}

\section{Further results on $\mathbb Z^n$}\label{Zn}\setcounter{equation}{0}

We now consider the $n-$dimensional {\it integer lattice graph}, i.e. $G=\mathbb Z^n$. We recall that, $x\sim y$ if and only if there exists $k\in \{1, \ldots, n\}$ such that $x_k=y_k\pm 1$ and $x_i=y_i$ for $i\neq k$. We define the edge weight and the node measure as
\[\omega: \mathbb Z^n\times \mathbb Z^n\to [0, +\infty);\quad\quad \omega(x,y)=\begin{cases}
1 & \text{ if } y\sim x\\
0 & \text{ if } y\not\sim x\,,
\end{cases}\]
\[\mu(x)=\sum_{y\in \mathbb Z^n} \omega(x,y)=2n\,.\]
We equip the graph $(\mathbb Z^n, \omega, \mu)$ with the euclidean distance
\begin{equation}\label{euc}
|x-y|=\left(\sum_{k=1}^n |x_k-y_k|^2 \right)^{\frac 12}\quad (x,y\in \mathbb Z^n)\,.\end{equation}

\begin{remark}
Observe that $\mathbb Z^n$ with the euclidean distance is not a weakly symmetric graph. In fact, in the definition of weakly symmetric graphs, only the combinatorial graph distance is considered. It is also easily seen that, $\mathbb Z^n$ endowed with the combinatorial metric, is not a weakly symmetric graph.
\end{remark}

On $\mathbb Z^n$, the condition on $\alpha$ made in \eqref{e12f} is not optimal.
In fact, the critical value is now $\alpha=2$, and no longer $\alpha=1$, as it will be clear from the next subsection.

\subsection{Phragm\`en-Lindel\"of principle and uniqueness}
In this case the condition on $\rho$ made in \eqref{e12f} (or more generally in \eqref{eq:rho-gen}) is not optimal. It turns out that it is indeed possible to consider even faster decaying densities. Let us set $x_0=0$, then we write $|x-x_0|=|x|$, i.e. the euclidean distance between $x$ and the reference point $x_0$. Here we assume that, for some  $\rho_0>0$ and $\alpha\in[0,2]$
\begin{equation}\label{eq:rho-Z}
\rho(x)\ge\rho_0(1+|x|)^{-\alpha}\,\,\,\text{for all}\,\,\,x\in \mathbb Z^n\,.
\end{equation}
More precisely, we can prove the next results.
\begin{theorem}\label{teo3}
Let $G=\mathbb Z^n$. Let $u$ be a subsolution of equation \eqref{problema} with $f\equiv0$, $u_0\equiv 0$ and $\rho$ such that \eqref{eq:rho-Z} holds. Furthermore, assume that $u$ fulfills \eqref{e3p}, with $x_0=0$, $d(x,y)$ being the euclidean distance \eqref{euc} and, for some $B>0$
\begin{equation}\label{eq56}
\overline Z(x):=
\begin{cases}
e^{B|x|^{\min\{1, 2-\alpha\}}} & \text{ if } \alpha\in [0, 2) \normalcolor \\
 e^{B\log^2(2+|x|^2)} & \text{ if } \alpha =2
\end{cases}, \text{ whenever } |x|\ge1.
\end{equation}
Then
\[u(x, t) \leq 0 \quad \forall (x,t)\in G_0^T.\]
\end{theorem}

A direct consequence of Theorem \ref{teo3} is the following uniqueness result.
\begin{corollary}\label{cor5}
Let $G=\mathbb Z^n$, $f\in \mathfrak F^T$ and $u_0\in\mathfrak F$. Assume that \eqref{eq:rho-Z} holds. Then there exists at most one solution to equation \eqref{problema} such that
\[\lim_{|x|\to +\infty}\frac1{Z(x)}\left\{\max_{t\in[0,T]}{|u(x,t)|}\right\}= 0\,,\]
where $Z$ is given by \eqref{eq56}.
\end{corollary}

\subsection{Optimality and nonuniqueness}
\begin{corollary}\label{cor3pa}
Let $G=\mathbb Z^n, n\geq 3$, and let
$u_0\in\mathfrak{F}$ be a \emph{non-negative and bounded} function on $G$.  Assume that
\begin{equation*}
\rho\in \mathfrak F, \quad 0<\rho(x)\leq c_0\, (1+|x|)^{-\alpha} \quad \textrm{for all}\;\; x\in G,
\end{equation*}
for some $\alpha>2$. Then for every fixed $\gamma\in \R$ satisfying
\eqref{eq:assumptionOnuzero},
there exists a solution $u$ to problem \eqref{problema} satisfying \eqref{eq:limituNonUnique}.
\end{corollary}

\section{Auxiliary Results}\label{auxiliary}\setcounter{equation}{0}

We now establish the following \emph{Weak Maximum Principle}.
\begin{lemma} \label{lem:WMP} Let assumption \eqref{e7f} be fulfilled.
 Let $\Omega\subseteq G$ be a finite set, and let $u\in\mathfrak{R}^T_0$ be such that
 \begin{equation} \label{eq:PBWMP}
  \begin{cases}
   \LL u \leq 0 & \text{in $\Omega^T$} \\
   u \leq 0 & \text{in $(G\setminus\Omega)^T_0$}\\
   u\leq0 &\text{in $\Omega\times\{0\}$}.
  \end{cases}
 \end{equation}
 Then 
 \begin{equation} \label{eq:toproveWMP}
	\text{$u\leq 0$ in $G^T_0$}. 
	\end{equation}
\end{lemma}
\begin{proof}
First of all, we observe that the assumption \(u\leq 0\) in
\((G\setminus\Omega)^T_0\) reduces the demonstration of
\eqref{eq:toproveWMP} to showing that
$
u\leq 0
$
in $\Omega^T_0$.
To this end, we argue essentially as in the proof of \cite[Lemma 3.3]{BiMePu}
(see also 
\cite[Lemma 1.39]{Grig2}).
We set $$M := \max_{\Omega\times[0,T]} u.$$ Observe that $M$ is well-defined since the set $\Omega\subseteq G$ is finite and $u\in\mathfrak{R}^T_0$ (hence, in particular,
 $u(x,\cdot)\in C([0,T])$ for all $x\in \Omega$). Then, let $(x_0,t_0)\in \Omega^T_0$ be such that
 $$u(x_0,t_0)=M.$$ 
 If $t_0 = 0$, then \eqref{eq:toproveWMP}
 follows from the fact that
 $u\leq 0$ on $\Omega\times\{0\}$. If, instead,  $t_0>0$, we assume by contradiction, that $M > 0$. Then, recalling that 
 $$\text{$\omega(x,y) > 0$ if $y\sim x$,}$$ 
 and taking into account \eqref{eq:PBWMP}, we have
 \begin{align*}
  0 \geq \LL u(x_0,t_0) &=\rho(x_0)\partial_tu(x_0,t_0)- \Delta u(x_0,t_0) \\
  &\ge -\frac{1}{\mu(x_0)}\sum_{y\in G}\omega(x_0,y)[u(y,t_0)-u(x_0,t_0)] \\
  & = \mathrm{Deg}(x_0)u(x_0,t_0)-\frac{1}{\mu(x_0)}\sum_{y\sim x_0}\omega(x_0,y)u(y,t_0)\\
  & = M\,\mathrm{Deg}(x_0)-\frac{1}{\mu(x_0)}\sum_{y\sim x_0}\omega(x_0,y)u(y,t_0)\,.
 \end{align*}
 Therefore, since $u\leq M$ in $\Omega^T_0$ and $u\leq 0 < M$ in $(G\setminus\Omega)^T_0$, we obtain
 $$M\,\mathrm{Deg}(x_0) \leq \frac{1}{\mu(x_0)} \sum_{y\sim x_0}\omega(x_0,y)u(y,t_0) \leq
  M\,\mathrm{Deg}(x_0),$$
 from which we derive that
 \begin{equation} \label{eq:ineqeqWMP}
 \frac{1}{\mu(x_0)}\sum_{y\sim x_0}\omega(x_0,y)u(y,t_0) = M\mathrm{Deg}(x_0).
  \end{equation}
  In view of \eqref{eq:ineqeqWMP}, and since $u\leq M$ in $G^T_0$, by
  \eqref{eq:DefDegWMP} we conclude that
  \begin{equation}\label{eq:connectedSWMP}
  u(y,t_0) = M \quad \text{for every}\,\,\, y\in G,\,\,\,y\sim x_0\,\,\,\text{with}\,\,\,u(x_0,t_0)=M.
  \end{equation}
Define $$F := \left\{(x,t_0),\,\, x\in G\,\,:\,u(x,t_0) = M\right\}.$$
 Now, let us consider some $(x,t_0)\in F$ and $y\in G\setminus \Omega$, hence $u(x,t_0)=M>0$ and $u(y,t_0)\le 0$. Due to \eqref{e7f}, there exist a path $\{x_k\}_{k=0}^n$ such that
 $$x_0 = x,\,\,\,x_n = y.$$ Since $x_0 = x$ and $(x,t_0)\in F$, we can apply \eqref{eq:connectedSWMP} and infer that $(x_1,t_0)\in F$. By repeating this argument, we get that $(x_i,t)\in F$ for every $i=0,...,n$; hence in particular that 
 $$(x_n,t_0) = (y,t_0)\in F$$ and thus $u(y,t_0)=M>0$ which yields a contradiction. This ends the proof.
 \end{proof}
\begin{remark}\label{rem:zeroImmaterial}
The choice of the initial time \(0\) in Lemma \ref{lem:WMP} \emph{is immaterial}.

More precisely, let \(0\le t_1<t_2<+\infty\), let \(\Omega\subseteq G\) be
finite, and let \(u\) have the same regularity as in Lemma \ref{lem:WMP} on
the time interval \([t_1,t_2]\). If
$$
\begin{cases}
\LL u\le0 & \text{in }\Omega\times(t_1,t_2] \\
u\le0 & \text{in }(G\setminus\Omega)\times[t_1,t_2] \\
u(\cdot,t_1)\le0 & \text{in }\Omega
\end{cases}
$$
then we can conclude that
\[
u\le0 \qquad\text{in }G\times[t_1,t_2].
\]
Indeed, it is enough to apply Lemma \ref{lem:WMP} to the time-translated
function
\[
\widetilde u(x,s):=u(x,t_1+s),
\qquad (x,s)\in G\times[0,t_2-t_1].
\]
\end{remark}
We now state a lemma which is used in Remark \ref{paragone-ell} and Section \ref{antitree}.

\begin{lemma}\label{lemmaelliptic}
Let there exists a function $Z\in \mathfrak F$, such that 
\begin{equation}\label{el1}
\Delta Z(x) \leq \rho(x) \quad \text{ for any } x\in G,
\end{equation} 
and, for some $c_0>0$
\begin{equation}\label{el3}
Z(x)\ge c_0\quad  \text{ for any } x\in G\,.
\end{equation}
Then, for $\gamma>\frac1{c_0}$,
$$\mathcal Z(x,t):=e^{\gamma t} Z(x),\quad (x,t)\in G^T_0$$
fulfills \eqref{e1p}. 
\end{lemma}

\begin{proof}
By \eqref{el3} we get
$$\mathcal Z(x,t)\geq {Z}(x) \ge c_0 > 0\quad\text{
for all $x\in G$ and $0<t\leq T$}.$$
This, together with \eqref{el1}, gives
\begin{equation*}
\begin{aligned}
\rho(x) \partial_t \mathcal Z - \Delta \mathcal Z & \geq  \rho(x) \gamma \mathcal Z(x,t) - \rho(x) \geq \rho(x)[\gamma{Z}(x) -1]\\ &\geq \rho(x)(\gamma\,c_0 -1)> 0 \quad \text{ for any } x\in G\,\,\,0<t\leq T,
\end{aligned}
\end{equation*}
provided that $\gamma > \frac1{c_0}$. This completes the proof.
\end{proof}

\section{Proofs of Proposition \ref{prop1}, Theorem \ref{teo1} and Theorem \ref{teo2}}\label{proofPL}\setcounter{equation}{0}

\begin{proof}[Proof of Proposition \ref{prop1}]
Let $\e > 0$ be arbitrarily fixed. 
On account of  \eqref{e3p}, we can infer the existence
of some $R_0 = R_0(\e)>0$ such that, for all $x$ with $d(x,x_0)>R_0$,
\begin{equation}\label{eq41}
\max_{t\in[0,T]}\frac{u(x,t)}{Z(x,t)}\,<\varepsilon\,.
\end{equation}
Then, we define
\[\mathcal Z_\varepsilon:= \varepsilon Z\,.\]
By assumption, it follows that, for any
fixed $R>R_0$, the function $\mathcal Z_\varepsilon$ is a supersolution of 
\begin{equation}\label{e12p}
\begin{cases}
\rho\partial_t u-\Delta u = 0 & \text{ in } B_{R}(x_0)\times(0,T]\\
u = \mathcal Z_\varepsilon & \text{ in } (G\setminus B_{R}(x_0))\times[0,T]\\
u = \mathcal Z_\varepsilon(\cdot,0) & \text{ in } B_{R}(x_0)\times\{0\}\,.
\end{cases}
\end{equation}
In fact, $\mathcal{Z}_\e\in\mathfrak{R}^T_0$ (as the same is true of $Z$), and by \eqref{e1p} we have
$$
\rho\,\partial_t \mathcal Z_\varepsilon-\Delta \mathcal Z_\varepsilon=\varepsilon\left(\rho\,\partial_t Z-\Delta Z\right)\ge0\qquad
\forall\,\,(x,t)\in B_{R}(x_0)\times(0,T].$$
Moreover, the exterior and initial inequalities required in Definition
\ref{defsolOmega} (see iii)-iv)) are identities, by the very definition of problem
\eqref{e12p}. 

On the other hand, $u$ is a subsolution of problem \eqref{e12p}. In fact, by assumption, $u$ is a sub\-so\-lution
of \eqref{problema}; thus, according to Definition \ref{defsol}, it satisfies
\begin{align*}
\ast)\,\,& u\in\mathfrak{R}^T_0\subseteq\mathfrak{R}^T_0(B_R(x_0)); \\
\ast)\,\,& \rho\partial_t u-\Delta  u\le0\,\,\,\text{in}\,\,\,G^T\supseteq B_R(x_0)\times(0,T]; \\
\ast)\,\,& u\le0\,(\le\mathcal Z_\varepsilon(\cdot,0))\,\,\,\text{in}\,\,G\times\{0\}\supseteq B_R(x_0)\times\{0\};
\end{align*}
where we have used the fact that $\mathcal Z_\varepsilon=\varepsilon Z>0$. 
Furthermore, due to \eqref{eq41}, we also have
$$
\frac{u(x,t)}{Z(x,t)}\le\max_{\tau \in[0,T]}\frac{u(x,\tau)}{Z(x,\tau)}<\varepsilon\qquad
\forall\,\,(x,t)\in (G\setminus B_{R}(x_0))\times[0,T]
$$
and therefore
$$
u(x,t)\le \varepsilon\,Z(x,t)=\mathcal Z_\varepsilon\quad \text{in}\,\,(G\setminus B_{R}(x_0))\times[0,T]\,.
$$
Hence, $u$ is a subsolution of \eqref{e12p}, as claimed. 
We can then apply Lemma \ref{lem:WMP} (to $u-\mathcal{Z}_\e$, which satisfies \eqref{eq:PBWMP} with $\Omega = B_R(x_0)$),
obtaining

$$u\leq \mathcal Z_ \varepsilon \quad \text{ in } G^T_0\,.$$
From this, by the arbitrariness of $\e > 0$, we conclude that
$$\text{$u\leq 0$ in $G^T_0$},$$
and the proof is complete. 
\end{proof}

The following Lemma, which will be useful in the proof of Theorems \ref{teo1} and \ref{teo2}, can be found in \cite[Lemma 5.1]{BP25}. We recall that $r$ has been defined in \eqref{def:r}.
\begin{lemma}
Let assumption \eqref{e7f} be satisfied. Let $\Omega\subset G$ be a finite set and let $f\in \mathfrak F$ be a spherically symmetric function with respect to $\Omega$.
Then
\begin{equation}\label{e13p}
\Delta f(x)=\mathfrak D_+(x)[f(r+1)-f(r)]+\mathfrak D_-(x)[f(r-1)-f(r)]
\end{equation}
for any $x\in G$ with $r\equiv r(x) \geq 1.$
\end{lemma}

\begin{proof}[Proof of Theorem \ref{teo1}]
Let $A > B$ be fixed, and let $Q >0$. We define
\begin{equation*}
\begin{gathered}
Z(x,t)= Z(r,t) := e^{A(1+Qt)(r+1)} \\[0.1cm]
(\text{with $r = r(x) = \bar{d}(x,\Omega)$})	
\end{gathered}
\end{equation*}
Clearly, \(Z\in\mathfrak{R}^{1/Q}_0\). We claim that one can choose
\(Q\) in such a way that \(Z\) fulfills the assum\-ptions of
Proposition \ref{prop1} for some \(x_0\in G\), \(d=\bar d\) and with $T$ replaced by
\begin{equation} \label{eq:DefTstar}
	T^* = \min\{1/Q,T\}.
\end{equation}
\noindent-\,\,\emph{Validity of \eqref{e1p}}.
In view of \eqref{e13p}, for all \((x,t)\in G^{1/Q}\) with \(r(x)\geq 1\),
\[
\begin{aligned}
\Delta Z(x,t)
&= \mathfrak D_+(x)\left[e^{A(1+Qt)(r+2)}-e^{A(1+Qt)(r+1)}\right] \\
&\quad - \mathfrak D_-(x)\left[e^{A(1+Qt)(r+1)}-e^{A(1+Qt)r}\right] \\
&=  \mathfrak D_+(x) Z(r,t)\left[e^{A(1+Qt)}-1\right]
 - \mathfrak D_-(x) Z(r,t)\left[1-e^{-A(1+Qt)}\right] \\
&\le \mathfrak D_+(x) Z(r,t)\left[e^{A(1+Qt)}-1\right].
\end{aligned}
\]
Therefore, for every \((x,t)\in G^{1/Q}\) with \(r(x)\geq 1\), using
\eqref{e12f}, we obtain
\[
\begin{aligned}
\rho\,\partial_t Z(x,t)- \Delta Z(x,t)
&\ge Z(r,t)\left\{\rho AQ(r+1)
-\mathfrak D_+(x)\left[e^{A(1+Qt)}-1\right]\right\}\\
&\ge Z(r,t)\mathfrak D_+(x)
\left\{\rho_0 AQ-\left[e^{A(1+Qt)}-1\right]\right\}\\
&\ge Z(r,t)\mathfrak D_+(x)
\left\{\rho_0 AQ-\left[e^{2A}-1\right]\right\} \geq 0,
\end{aligned}
\]
for every $x\in G$ with $r = r(x)\geq 1$ and every $0\leq t\leq 1/Q$, provided that
\[
Q\ge\frac{e^{2A}-1}{\rho_0A},
\]
In order to extend this estimate to the whole of $G$ (hence, to all $x\in \Omega$) it suffices to observe that, since 
\emph{$\Omega$ is finite and $G$ is locally finite}, we have
$$0\leq r(y) \leq \mathfrak{r}\quad\text{for every
$x\in \Omega$ and every $y\sim x$},$$
for some $\mathfrak{r} > 0$ only depending on $\Omega$.
Thus, we obtain
\begin{align*}
\Delta Z(x,t) & = \frac{1}{\mu(x)}\sum_{y\in G}
[Z(r(y),t)-Z(r(x),t)]\omega(x,y) \\
& \leq 	e^{2A(\mathfrak{r}+1)}\mathrm{Deg}(x)
\leq e^{2A(\mathfrak{r}+1)}\cdot\max_{x\in\Omega}\mathrm{Deg}(x) \equiv \kappa.
\end{align*}
From this, we conclude that
\begin{equation} \label{eq:DaRichiamareZn}
\begin{split}
	\rho(x)\partial_t Z(x,t)-\Delta Z(x,t) 
& = \rho(x)AQ e^{A(1+Qt)}-\Delta Z(x,t) \\
& \geq \rho(x)AQe^A-\kappa \\
& \geq Q\cdot\big(Ae^A\min_{x\in\Omega}\rho(x)\big)-\kappa \geq 0,
\end{split}	
\end{equation}
for every $x\in \Omega$ (hence, with $r= r(x) = 0$) and every $0\leq t\leq 1/Q$, provided that 
$$Q \geq Q(A) = \kappa\cdot\big(Ae^A\min_{x\in\Omega}\rho(x)\big)^{-1}. $$
(recall that $\rho > 0$ on $G$). Summing up, we have
$$\rho(x)\partial_t Z(x,t)-\Delta Z(x,t)\geq 0 \quad \text{on $G^{1/Q}$},$$
provided that $Q$ is large enough, and thus
$Z$ satisfies \eqref{e1p}, as claimed. We stress that this lower bound on \(Q\)
is independent of the final time \(T\); this will allow us to iterate the
argument below on time intervals of length at most \(1/Q\).
\medskip

\noindent-\,\,\emph{Validity of \eqref{e3p}}.
Let \(T^*\) be as in \eqref{eq:DefTstar}, with \(Q\) chosen as above.
We first observe that, since $u$ satisfies \eqref{e16p}
(and by
Remark \ref{rem:Equivalenteuplus}), for every $\e > 0$
there exists $r_0 > 0$ such that
$$\frac{1}{\tilde{Z}(x)}\max_{t\in[0,T]}u^+(x,t) <\e \quad
\text{for all $x\in G$ with $r = r(x)\geq r_0$}.$$
On the other hand, if we arbitrarily fix $x_0\in \Omega$
and if set $d_0 = \mathrm{diam}(\Omega)<+\infty$
(recall that $\Omega$ is finite), for every $x\in G$ with $d(x,x_0) > r_0+d_0$ and every $y\in\Omega$ we have
$$d(x,y)\geq d(x,x_0)-d(x_0,y)\geq r_0,$$
and thus $r = r(x)\geq r_0$; as a consequence, we get
\begin{equation} \label{eq:LimituplusdFirst}
\frac{1}{\tilde{Z}(x)}\max_{t\in[0,T]}u^+(x,t) <\e \quad
\text{for all $x\in G$ with $d(x,x_0)>r_0+d_0$}.	
\end{equation}
Using \eqref{eq:LimituplusdFirst}, and since $A > B$,
we get
\begin{align*}
\max_{t\in [0,T^*]}\frac{u^+(x,t)}{Z(r,t)}
& \leq \max_{t\in [0,T]}\frac{u^+(x,t)}{Z(r,t)}
\leq \frac{1}{e^{A(r+1)}}\max_{t\in[0,T]}u^+(x,t)	 \\
& = e^{(B-A)(r+1)}\Big(\frac{1}{\tilde{Z}(x)}\max_{t\in[0,T]}u^+(x,t)\Big) < \e,
\end{align*}
for every $x\in G$ with $d(x,x_0) > r_0+d_0$, 
and hence \(Z\) satisfies \eqref{e3p} with $T$ replaced by
$T^*$
(actually, $Z$ satisfies the condition
in Remark \ref{rem:Equivalenteuplus}, which is equivalent
to \eqref{e16p}).
\vspace{0.1cm}

Gathering all these facts, we can exploit
Proposition \ref{prop1}, obtaining
$
u\le0
$
in $G_0^{T^*}$.
By ite\-rating this argument, after finitely many steps 
(since \(T<+\infty\) and since \(Q\), being independent of \(T\), can be kept fixed
throughout the iteration) we finally get
\[u\le0
\qquad\text{in }G_0^T.
\]
This is the desired conclusion.
\end{proof}

\begin{proof}[Proof of Theorem \ref{teo2}]  Let 
$A\in (B,\rho_0/2)$ be fixed (here $\rho_0 > 0$
is as in assumption \eqref{eq:rho-gen},
and $B<\rho_0/2$ is as in \eqref{eq:soprasol}), 
and let $Q >0$. We define
\begin{equation*}
\begin{gathered}	
	Z(x,t) = Z(r,t):=\exp\left(A(1+Qt)(r+1)\log^\beta(r+2)\right),
	\\[0.1cm]
	\text{(where $r = r(x) = \bar{d}(x,\Omega)$)}.
\end{gathered}
\end{equation*}
Clearly, \(Z\in\mathfrak{R}^{1/Q}_0\). We claim that one can choose
\(Q\) in such a way that \(Z\) fulfills the assum\-ptions of
Proposition \ref{prop1} for some \(x_0\in G\),  \(d=\bar d\)
and with $T$ replaced by
\begin{equation} \label{eq:DefTstarBIS}
	T^* = \min\{1/Q,T\}.
\end{equation}
Throughout what follows, to simplify the notation we set
\begin{align*}
i)\,\,&\Phi(s):=(s+1)\log^\beta(s+2) & &(s\ge0);	\\
ii)\,\,&k(t):=A(1+Qt) & & (0\leq t\leq 1/Q).
\end{align*}
\noindent -\,\,\emph{Validity of \eqref{e1p}}.
In view of \eqref{e13p}, and since the function
$\Phi$ is \emph{smooth and strictly increasing}
on the interval $[0,+\infty)$, for every $x\in G$ with $r = r(x)\geq 1$ we have
\begin{align*}
	\Delta Z(x,t)
&=
\mathfrak D_+(x)\bigl[Z(r+1,t)-Z(r,t)\bigr]
-
\mathfrak D_-(x)\bigl[Z(r,t)-Z(r-1,t)\bigr]  \\
&\le
\mathfrak D_+(x)\bigl[Z(r+1,t)-Z(r,t)\bigr] \\
& \text{(by the Mean Value theorem)} \\
& = 
\mathfrak{D}_+(x)\big[Z(\eta,t)
k(t)\Phi'(\eta)\big],
\end{align*}
for some \(\eta\in[r,r+1]\).
We then observe that, by definition,
$$\Phi'(s) = \log^\beta(s+2)+\beta\cdot\frac{s+1}{s+2}\log^{\beta-1}(s+2).$$
We stress that estimating \(\Phi'\) in terms of
\(\log^\beta(s+1)\), rather than \(\log^\beta(s+2)\), is useful in what
follows, since we shall apply the Mean Value Theorem with
\(s=\eta\in[r,r+1]\), and hence
\[
\eta+1\le r+2.
\]
As a consequence, since we clearly have $\Phi'(s)\sim \log^\beta(s+1)$ as $s\to+\infty$,
for every $\sigma > 0$ there exists $R = R_\sigma > 0$ such that  
\begin{equation} \label{eq:estimPhiprimetouse}
\Phi'(s)\leq (1+\sigma)\log^\beta(s+1)\qquad\forall\,\,s\geq R_\sigma.	
\end{equation}
In particular, from \eqref{eq:estimPhiprimetouse} we also derive
\begin{equation} \label{eq:MeanValuePhi}
	0\leq \Phi(s_2)-\Phi(s_1) \leq (1+\sigma)\log^\beta(s_2+1)(s_2-s_1)\quad\forall\,\,s_2\geq s_1\geq R_\sigma.
\end{equation}
Now, choosing  $0<\sigma<1$ so small that
$$2A(1+\sigma)<\rho_0
$$
(notice that this choice is possible, since $A<\rho_0/2$),
we obtain
\begin{equation} \label{eq:stimaDeltaZ}
\begin{split}
\Delta Z(x,t) & \leq 2\mathfrak{D}_+(x)k(t)Z(\eta,t)\log^\beta(\eta+1) \\
& (\text{by definition of $Z$}) \\
& = 2\mathfrak{D}_+(x)k(t)Z(r,t)e^{k(t)(\Phi(\eta)-\Phi(r))}\log^\beta(\eta+1) \\
& (\text{since $r<\eta<r+1$, and using \eqref{eq:MeanValuePhi}}) \\
& \leq 2\mathfrak{D}_+(x)k(t)Z(r,t)e^{k(t)(1+\sigma)\log^\beta(r+2)}\log^\beta(r+2) \\
& \leq 4A\mathfrak{D}_+(x)Z(r,t)e^{2A(1+\sigma)\log^\beta(r+2)}\log^\beta(r+2) \\
& \leq 4A\mathfrak{D}_+(x)Z(r,t)e^{\rho_0\log^\beta(r+2)}\log^\beta(r+2),
\end{split}	
\end{equation}
and this estimate holds for every $x\in G$ with $r = r(x)\geq R_\sigma$, and every $0\leq t\leq 1/Q$.

Thanks to \eqref{eq:stimaDeltaZ}, and using
assumption \eqref{eq:rho-gen}, we finally get
\begin{align*}
& \rho(x)\partial_t Z(x,t)-\Delta Z(x,t) 
= \rho(x)AQ\Phi(r)Z(r,t)-\Delta Z(r,t) \\
& \quad \geq Z(r,t)\log^\beta(r+2)\big[\rho(x)AQ(r+1)-4A\mathfrak{D}_+(x)e^{\rho_0\log^\beta(r+2)}\big] \\
& \quad \geq \mathfrak{D}_+(x)Z(r,t)\log^\beta(r+2)e^{\rho_0\log^\beta(r+2)}\big[A(Q-4)\big]
\geq 0, 
\end{align*}
for every $x\in G$ with $r = r(x)\geq R_\sigma$ and every
$0\leq t\leq 1/Q$, provided that
$$Q \geq 4.$$
In order to extend this estimate to the whole of $G$ (hence, to all $x\in G$ with $r = r(x)\leq R_\sigma$) it suffices to observe that, since 
\emph{$\Omega$ is finite and $G$ is locally finite}, we have
$$0\leq r(y) \leq \mathfrak{r}\quad\text{for every
$x\in G$ with $r(x)\leq R_\sigma$ and every $y\sim x$},$$
for some $\mathfrak{r} > 0$ only depending on $\Omega$ and $R_\sigma$
(see also Remark \ref{rem:SferaOmegafinita}).

Thus, we obtain
\begin{align*}
\Delta Z(x,t) & = \frac{1}{\mu(x)}\sum_{y\in G}
[Z(r(y),t)-Z(r(x),t)]\omega(x,y) \\
& \leq 	e^{2A\Phi(\mathfrak{r})}\mathrm{Deg}(x)
\leq e^{2A\Phi(\mathfrak{r})}\cdot\max_{r(x)\leq R_\sigma}\mathrm{Deg}(x) \equiv \kappa.
\end{align*}
From this, we conclude that
\begin{align*}
	\rho(x)\partial_t Z(x,t)-\Delta Z(x,t) 
& = \rho(x)AQ\Phi(r)Z(r,t)-\Delta Z(r,t) \\
& \geq \rho(x)QA\Phi(r)Z(r,t)-\kappa \\
& \geq Q\cdot\big(A\Phi(0)e^{A\Phi(0)}\min_{r(x)\leq R_\sigma}\rho(x)\big)-\kappa \geq 0,
\end{align*}
for every $x\in G$ with $r = r(x)\leq R_\sigma$ and every $0\leq t\leq 1/Q$, provided that 
$$Q \geq Q(A;R_\sigma) = \kappa\cdot\big(A\Phi(0)e^{A\Phi(0)}\min_{r(x)\leq R_\sigma}\rho(x)\big)^{-1}. $$
(recall that $\rho > 0$ on $G$). Summing up, we have
$$\rho(x)\partial_t Z(x,t)-\Delta Z(x,t) \geq 0\quad \text{on $G^{1/Q}$},$$
provided that $Q$ is large enough,
and thus
$Z$ satisfies \eqref{e1p}, as claimed.
We stress that this lower bound on \(Q\)
is independent of the final time \(T\); this will allow us to iterate the
argument below on time intervals of length at most \(1/Q\).
\medskip

\noindent-\,\,\emph{Validity of \eqref{e3p}}.
Let \(T^*\) be as in \eqref{eq:DefTstarBIS}, with \(Q\) chosen as above.
We first observe that, since $u$ satisfies \eqref{eq:limite}
(and by
Remark \ref{rem:Equivalenteuplus}), for every $\e > 0$
there exists $r_0 > 0$ such that
$$\frac{1}{\hat{Z}(x)}\max_{t\in[0,T]}u^+(x,t) <\e \quad
\text{for all $x\in G$ with $r = r(x)\geq r_0$}.$$
On the other hand, if we arbitrarily fix $x_0\in \Omega$
and if set $d_0 = \mathrm{diam}(\Omega)<+\infty$
(recall that $\Omega$ is finite), for every $x\in G$ with $d(x,x_0) > r_0+d_0$ and every $y\in\Omega$ we have
$$d(x,y)\geq d(x,x_0)-d(x_0,y)\geq r_0,$$
and thus $r = r(x)\geq r_0$; as a consequence, we get
\begin{equation} \label{eq:Limituplusd}
\frac{1}{\hat{Z}(x)}\max_{t\in[0,T]}u^+(x,t) <\e \quad
\text{for all $x\in G$ with $d(x,x_0)>r_0+d_0$}.	
\end{equation}
Using \eqref{eq:Limituplusd}, and since $A > B$
(and $Q\geq 1/T$), we then have
\begin{align*}
\max_{t\in [0,1/Q]}\frac{u^+(x,t)}{Z(r,t)}
& \leq \max_{t\in [0,T]}\frac{u^+(x,t)}{Z(r,t)}
\leq \frac{1}{e^{A(r+1)\log^\beta(r+2)}}\max_{t\in[0,T]}u^+(x,t)	 \\
& = e^{(B-A)(r+1)\log^\beta(r+2)}\Big(\frac{1}{\hat{Z}(x)}\max_{t\in[0,T]}u^+(x,t)\Big) < \e,
\end{align*}
for every $x\in G$ with $d(x,x_0) > r_0+d_0$, 
and thus \(Z\) satisfies \eqref{e3p}
with $T$ replaced by
$T^*$
(actually, $Z$ satisfies the condition
in Remark \ref{rem:Equivalenteuplus}, which is equivalent
to \eqref{eq:limite}).
\vspace{0.1cm}

Gathering all these facts, we can exploit
Proposition \ref{prop1}, obtaining
$
u\le0
$
in $G_0^{T^*}$.
By ite\-rating this argument, after finitely many steps 
(since \(T<+\infty\) and since \(Q\), being independent of \(T\), can be kept fixed
throughout the iteration) we finally get
\[
u\le0
\qquad\text{in }G_0^T.
\]
This is the desired conclusion.
\end{proof}
\section{Proofs of Theorem \ref{prop71}, Corollary \ref{cor1pa}}\label{nonuniqueness}\setcounter{equation}{0}
To prove Theorem \ref{prop71}, we first show the following existence result.
 \begin{proposition} \label{prop:SolvDir} Let assumption \eqref{e7f} be in force,
 and let $\rho\in\mathfrak{F}$ be such that $\rho > 0$ on $G$.
 Moreover, let $\Omega\subseteq G$ be a 
 \emph{finite set}, and let
 $T\in(0,+\infty)$.
 Finally, let $$u_0\in\mathfrak{F}(\Omega),$$ 
 and let $f,g$ be two functions satisfying the following properties:
 \begin{itemize}
  \item[i)] $f\in \mathfrak{F}^T_0(\Omega)$, and
  $f(x,\cdot)\in C([0,T])$ for all $x\in\Omega$;
  \vspace*{0.1cm}

  \item[ii)] $g\in\mathfrak{F}^T_0(G\setminus\Omega)$, 
  and 
  $g(x,\cdot)\in C([0,T])$ for all $x\notin\Omega$;
  
\end{itemize}
Then there exists a \emph{unique solution} $u\in \mathfrak{R}^T_0(\Omega)$ to problem
\eqref{eq:pbDir}, namely
\begin{equation*}
	\begin{cases}
\LL u = f & \text{in $\Omega^T$} \\
u = g & \text{in $(G\setminus\Omega)^T_0$} \\
u = u_{0} & \text{in $\Omega\times\{0\}$}
\end{cases}
\end{equation*}
(in the sense of Definition \ref{defsolOmega}). 
\end{proposition}
\begin{proof}
We begin by proving the \emph{uniqueness part} of the proposition. To this end, let us
 assume that there exist
 two solutions $u_1,u_2\in\mathfrak{R}^T_0(\Omega)$ of problem \eqref{eq:pbDir}
 (in the sense of Definition
 \ref{defsolOmega}), and let
 $w = u_1-u_2\in\mathfrak{R}^T_0(\Omega).$
 Since both $u_1$ and $u_2$ solve \eqref{eq:pbDir}, we clearly have
 \begin{itemize}
  \item[i)] $\LL w = f-f = 0$ on $\Omega^T$;
  \item[ii)] $w(x,t) = g(x,t)-g(x,t) = 0$ on $(G\setminus\Omega)^T_0$;
  \item[iii)] $w(x,0) = u_0(x)-u_0(x) = 0$ for all $x\in \Omega$.
 \end{itemize}
 As a consequence, by applying the Weak Maximum Principle in Lemma \ref{lem:WMP} to $\pm w$, we
 conclude that $w = 0$ pointwise on $G^T_0$, and therefore 
 $$\text{$u_1 = u_2$ on $G_0^T$}.$$
 \noindent We now turn to prove the
 \emph{existence part}.
 To begin with, we
 observe that, since
 $\Omega$ is finite, we can write $\Omega = \{x_1,\ldots,x_p\}$ (for some $p\in\mathbb{N}$);
 thus, given any $u\in\mathfrak{R}^T_0(\Omega)$,
  we see that
 $u$ is a solution to problem
 \eqref{eq:pbDir} \emph{if and only if}
 \begin{itemize}
 	\item[1)] $u(x,t) = g(x,t)$ for all $x\in G\setminus\Omega,\,t\in[0,T]$;
 	\item[2)] $u(x,0) = u_0(x)$ for all $x\in \Omega$; 
 	\item[3)] for every $1\leq i\leq p$ and every
 	$0<t\leq T$, one has
 	\begin{align*}
 	\partial_t u(x_i,t) & = \frac{1}{\rho(x_i)}\Delta u(x_i,t)+\frac{1}{\rho(x_i)}f(x_i,t)\\
 	& = \frac{1}{\rho(x_i)}\cdot\frac{1}{\mu(x_i)}\sum_{y\in G}
 	[u(y,t)-u(x_i,t)]\omega(x_i,y)+\frac{1}{\rho(x_i)}f(x_i,t) \\
 	& = -\frac{\mathrm{Deg}(x_i)}{\rho(x_i)}u(x_i,t)+
 	\frac{1}{\rho(x_i)}\cdot\frac{1}{\mu(x_i)}\sum_{y\in G}
 	u(y,t)\omega(x_i,y) +\frac{1}{\rho(x_i)}f(x_i,t)\\
 	& (\text{using the boundary condition in 1)}) \\
 	& = -\frac{\mathrm{Deg}(x_i)}{\rho(x_i)}u(x_i,t)+
 	\frac{1}{\rho(x_i)}\cdot\frac{1}{\mu(x_i)}\sum_{j=1}^p
 	u(x_j,t)\omega(x_i,x_j)\\
 	& \qquad\qquad+\frac{1}{\rho(x_i)}\cdot\frac{1}{\mu(x_i)}\sum_{y\notin\Omega}g(y,t)\omega(x_i,y)
 	+\frac{1}{\rho(x_i)}f(x_i,t).
  	\end{align*}
 \end{itemize}
 In other words, if $u\in\mathfrak{R}^T_0(\Omega)$
 satisfies the boundary condition in 1), the above identity shows that $u$ solves \eqref{eq:pbDir} 
 \emph{if and only if}
 the functions
 $$\mathfrak{u}_i(t):[0,T]\to\R,\qquad 
 \mathfrak{u}_i(t) = u(x_i,t)\qquad (1\leq i\leq p)$$
  satisfy the following \emph{linear 
  Cauchy problem in $\R^p$}:
\begin{equation} \label{eq:PbCauchyAusiliario}
	\begin{cases} 
z' = Az+ F(t) \\
z(0) = \mathfrak{u}_0	
\end{cases}
\end{equation}
where $A = (a_{i,j})\in M_p(\R),\,F = (F_1,\ldots,F_p):[0,T]\to\R^p$ and $\mathfrak{u}_0\in\R^p$ are given by
\begin{align*}
\ast)\,\,& a_{i,j} = 	-\frac{\mathrm{Deg}(x_i)}{\rho(x_i)}\delta_{i,j}+\frac{1}{\rho(x_i)}\cdot\frac{1}{\mu(x_i)}\omega(x_i,x_j)\qquad (1\leq i,j\leq p); \\[0.1cm]
\ast)\,\,& F_i(t) =  \frac{1}{\rho(x_i)}\cdot\frac{1}{\mu(x_i)}\sum_{y\notin\Omega}g(y,t)\omega(x_i,y)+\frac{1}{\rho(x_i)}f(x_i,t)\qquad (1\leq i\leq p); \\
\ast)\,\,&\mathfrak{u}_0 = (u_0(x_1),\ldots,u_0(x_p))
\end{align*}
(here, $\delta_{i,j}$ stands for the usual Kronecker symbol). 
\vspace{0.1cm}

Now, by the assumptions on $f,g$ (and since the graph $G$ is locally finite, see \eqref{e7f}), one can easily recognize
that 
$F\in C([0,T];\R^p)$; thus by classical results of ODE Theory we infer that
the Cauchy problem
\eqref{eq:PbCauchyAusiliario} possesses a unique solution 
$$\mathfrak{u} = (\mathfrak{u}_1,\ldots, \mathfrak{u}_p)\in C^1([0,T];\R^p).$$
Then, if we define 
$$u:G^T_0\to\R,\qquad u(x,t) = \begin{cases}
 g(x,t) & \text{if $x\notin\Omega$} \\
 \mathfrak{u}_i(t) & \text{if $x\in\Omega,\,x = x_i$}	
 \end{cases}
$$
we deduce that $u$ is a solution to problem
\eqref{eq:pbDir} 
(in the sense of Definition \ref{defsolOmega}).
\vspace{0.05cm}

 Indeed, since $\mathfrak{u}\in C^1([0,T];\R^p)$, it follows that
$u\in\mathfrak{R}^T_0(\Omega);$
moreover, by definition we also have 
$u=g$ on $(G\setminus\Omega)^T_0$
and $u= u_0$ on $\Omega\times\{0\}$.
Finally, since $\mathfrak{u}$
solves
\eqref{eq:PbCauchyAusiliario}, from the above computation 
we infer that
$\text{$\LL u = f$ in $\Omega^T$}.$
\vspace{0.05cm}

Gathering all these facts, we conclude that
$u$ solves \eqref{eq:pbDir}, as desired.
 \end{proof}
 
\begin{remark} \label{rem:HpThesisExistence}
Although Definition \ref{defsolOmega}, for solutions of problem \eqref{eq:pbDir},
only requires \emph{time regu\-la\-rity at the points of $\Omega$}, it follows
from the proof of Proposition \ref{prop:SolvDir}
\emph{(}of which we inherit the notation\emph{)} that the solution
$u\in\mathfrak{R}^T_0(\Omega)$ constructed there enjoys the following
properties:
\begin{itemize}
\item for every $x\in\Omega$, we have
$
u(x,\cdot)\in C^1([0,T])
$,
and $\LL u=f$
 in $\Omega^T_0$;
\vspace{0.1cm}

\item for every $x\in G\setminus\Omega$, we have
$
u(x,\cdot)\in C([0,T]).
$
\end{itemize}
The first property follows from the fact that the function $f$ is assumed
to be defined and conti\-nuous up to $t=0$, together with the continuity up to
$t=0$ of the function $g$, since both enter the forcing term $F(t)$ of the
finite-dimensional ODE system associated with the problem. 

The second property
follows directly from the boundary condition
\[
u=g\qquad\text{in }(G\setminus\Omega)^T_0
\]
and from the assumed continuity of $g$ on $[0,T]$.
\end{remark}

With the above results at hand, we can now prove Theorem \ref{prop71}.
\begin{proof}[Proof of Theorem \ref{prop71}] 
Let $B = B_{\hat{R}}(o)\subseteq G$ and $u_0\in\mathfrak{F}$
be as in the statement of the theorem. We arbitrarily fix $\gamma\in\R$ satisfying condition \eqref{eq:assumptionOnuzero}, namely
$$\gamma\geq M_0 = \sup_G u_0$$
(recall that $0\leq M_0<+\infty$, since $u_0$ is non-negative and bounded), and we 
use an \emph{exhaustion-plus-barrier} argument to
construct a solution $u_\gamma$
to  \eqref{problema} satisfying  \eqref{eq:limituNonUnique}, that is,
$$\text{$u_\gamma(x,t_0)\to \gamma$ as $d(x,o)\to+\infty$,\quad
for each fixed $t_0 > 0$}.$$
In order to ease the readability, we proceed by steps.
\medskip

\textsc{Step I): Construction of $u_\gamma$}. 
To begin with, for every $j\in\mathbb{N}$
 we consider
the following \emph{approximating} Cauchy-Dirichlet problem for $\LL$
\begin{equation}\label{eq711}
\begin{cases}
\LL u=0 & \text{in}\,\,B_j(o)\times(0,T]\\
u = \gamma & \text{in}\,\, (G\setminus B_j(o))\times[0,T] \\
u(x,0) = u_0(x) & \text{for every $x\in B_j(o)$}.
\end{cases}
\end{equation}
On account of Proposition \ref{prop:SolvDir} (applied here
with $f \equiv 0$ and $g\equiv \gamma$), we know that
this problem possesses a \emph{unique solution $v_j
\in\mathfrak{R}_0^T(B_j(o))$}, further satisfying the following properties:
\begin{itemize}
\item for every $x\in B_j(o)$, we have
$
v_j(x,\cdot)\in C^1([0,T])
$,
and $\LL v_j=0$
 in $B_j(o)\times[0,T]$;
\vspace{0.1cm}

\item for every $x\in G\setminus B_j(o)$, we have
$
v_j(x,\cdot)\in C([0,T]).
$
\end{itemize}
(see Remark \ref{rem:HpThesisExistence}).
We then claim that the following facts hold.
\begin{itemize}
\item[(1)] For every $(x,t)\in G_0^T$ and every $j\in\mathbb{N}$, one has
\begin{equation}\label{eq712}
0\le v_j(x,t)\le \gamma;
\end{equation}
\item[(2)] The sequence $\{v_j\}_j$ is \emph{decreasing}.
\end{itemize}

\noindent -\,\,\emph{Proof of Claim} (1). On the one hand,
since $v_j\in\mathfrak{R}_0^T(B_j(o))$ solves problem
\eqref{eq711}, and since both the (constant) exterior datum $\gamma$ and the initial datum $u_0$ are non-negative,
by the Weak Maximum Principle in Lemma \ref{lem:WMP} 
(applied here to $-v_j$)
we derive that
$$\text{$v_j(x,t) \geq 0$ for every $(x,t)\in G_0^T$}.$$
On the other hand, 
again recalling that $v_j$ solves \eqref{eq711}, 
and since $\gamma\geq M_0\geq u_0$ on $G$, 
by the Weak Maximum Principle in Lemma \ref{lem:WMP} 
(applied here to $v_j-\gamma$) we get
$$\text{$v_j(x,t)\leq \gamma$ for every $(x,t)\in G_0^T$}.$$
Hence, Claim (1) is proved.
\vspace*{0.1cm}

\noindent-\,\,\emph{Proof of Claim} (2). 
We apply once again the Weak Maximum Principle
in Lemma \ref{lem:WMP}. First of all, since
$v_j$ is a solution of problem \eqref{eq711}, setting $w_j = v_{j+1}-v_j$ we  have
$$
\LL w_j(x,t) = 0\quad\text{for every $(x,t)\in B_j(o)\times(0,T]$}.
$$
Moreover, on account of \eqref{eq712} we also get
\begin{itemize}
 \item $w_j(x,t) = v_{j+1}(x,t)-v_j(x,t) = v_{j+1}(x,t)-\gamma \leq 0$ on $(G\setminus B_j(o))\times[0,T]$;
 \vspace*{0.1cm}

 \item $w_j(x,0) = v_{j+1}(x,0)-v_j(x,0) = 0$ for every $x\in B_j(o)\subseteq B_{j+1}(o)$.
\end{itemize}
Therefore, by Lemma \ref{lem:WMP}, $w_j\leq 0$ in $G_0^T$; hence, for any $j\in\mathbb{N}$, we have
\begin{equation*}
v_{j+1}\leq  v_j\quad \text{in}\,\,\,G^T_0,
\end{equation*}
and this completes the proof of Claim (2).
\vspace*{0.1cm}

Now, by combining Claim (1) and Claim (2) we deduce that the sequence $\{v_j\}_{j}$ is \emph{decre\-asing and bounded} on $G_0^T$; therefore, there exists $u_\gamma\in\mathfrak{F}_0^T$ such that
\begin{equation} \label{eq:ugammaDefBd}
\begin{split}
  \ast)\,\,&\text{$u_\gamma(x,t) = \lim_{j\to+\infty}v_j(x,t)$ for every $(x,t)\in G_0^T$}; \\
   \ast)\,\,&\text{$0\leq u_\gamma(x,t)\leq \gamma$ for every $(x,t)\in G_0^T$}.	
\end{split}
\end{equation}
In particular, we also have
\begin{equation}\label{eq:DeltavjtoDeltaugamma}
	\text{$\Delta v_j(\bar{x},\cdot)\to \Delta u_\gamma(\bar{x},\cdot)$ pointwise on $[0,T]$, as $j\to+\infty$}.
\end{equation}
(since the sum
defining $\Delta$ is a \emph{finite sum}, see
\eqref{e7f} and recall that $\omega(x,y) > 0\Leftrightarrow x\sim y$).
\vspace{0.1cm}

\textsc{Step II): Regularity of $u_\gamma$.} We now
fix $\bar{x}\in G$, and we turn to prove that
\begin{equation} \label{eq:ToProveugammaC1}
	u_\gamma(\bar{x},\cdot)\in C^1([0,T]).
\end{equation}
To begin with, we show that $u_\gamma(\bar{x},\cdot)$ is
\emph{continuous} on $[0,T]$ by using a \emph{compactness argument} based on the Arzel\'a-Ascoli theorem. 
On the one hand, from \eqref{eq712} we immediately
derive that the sequence $\{v_j(\bar{x},\cdot)\}_j$
is \emph{equi-bounded on $[0,T]$}. On the other hand, 
since
$$\text{$v_j(\bar{x},\cdot)\in C^1([0,T])$ 
\quad and \quad $\LL v_j = 0$ on $B_j(o)\times[0,T]$}$$
(see \textsc{Step I)}), for every $j > d(\bar{x},o)$ and every
$0\leq s,t\leq T$ we have
\begin{equation} \label{eq:identityToUseForEquation}
\begin{split}
	 v_j(\bar{x},t)-v_j(\bar{x},s) & = \int_s^t	\partial_tv_j(\bar{x},\tau)\,d\tau  = \frac{1}{\rho(\bar{x})}\int_s^t \Delta v_j(\bar{x},\tau)\,d\tau,\end{split}
\end{equation}
and $\Delta v_j(\bar{x},\cdot)$ is \emph{bounded on $[0,T]$}:
indeed, using once again \eqref{eq712}, we get
\begin{equation} \label{eq:ToPassLimitDCT}
\begin{split}
	|\Delta v_j(\bar{x},\tau)| & \leq 
	\frac{1}{\mu(\bar{x})}
	\sum_{y\in G}|v_j(y,\tau)-v_j(\bar{x},\tau)|\omega(\bar{x},y) \\
	& (\text{recall \eqref{eq:DefDegWMP}}) \\
	& \leq 2\gamma\mathrm{Deg}(\bar{x})\quad\forall\,\,0\leq \tau\leq T.
\end{split}
\end{equation}
As a consequence, we obtain
$$|v_j(\bar{x},t)-v_j(\bar{x},s)|\leq \frac{2\gamma\mathrm{Deg}(\bar{x})}{\rho(\bar{x})}\cdot|t-s|,$$
and this proves that $\{v_j(\bar{x},\cdot)\}_j$ is \emph{equicontinuous on $[0,T]$}.
\vspace{0.1cm}

Gathering all these facts, and since we \emph{already know that $v_j(\bar{x},\cdot)\to u_\gamma(\bar{x},\cdot)$} pointwise on (the compact interval) $[0,T]$ as $j\to+\infty$, by
the Arzel\'a-Ascoli Theorem we get
\begin{equation} \label{eq:contugamma}
u_\gamma(\bar{x},\cdot)\in C([0,T]).	
\end{equation}
In particular, from \eqref{eq:ugammaDefBd} we also derive
that, for every $x\in G$,  
\begin{equation} \label{eq:ToUseEquationugamma}
\Delta u_\gamma({x},\cdot) = 
\frac{1}{\mu({x})}
	\sum_{y\in G}[u_\gamma(y,\cdot)-u_\gamma({x},\cdot)]\omega({x},y) 
\in C([0,T])	
\end{equation}
(as the sum defining $\Delta$ is a \emph{finite sum}, and $\bar{x}$ in \eqref{eq:contugamma}
 is arbitrary).
\vspace{0.1cm}

Now, starting from identity \eqref{eq:identityToUseForEquation},
and using \eqref{eq:DeltavjtoDeltaugamma}, \eqref{eq:ToPassLimitDCT} and \eqref{eq:ToUseEquationugamma}, we can easily complete the proof of \eqref{eq:ToProveugammaC1}.
Indeed, using the cited \eqref{eq:identityToUseForEquation} with $s = 0$ (and recalling that $v_j$ is a solution of problem \eqref{eq711}), 
for every $j > d(\bar{x},o)$ and every $0\leq t\leq T$ we can write 
$$v_j(\bar{x},t) = u_0(\bar{x}) +\frac{1}{\rho(\bar{x})}\int_0^t \Delta v_j(\bar{x},\tau)\,d\tau;$$
as a consequence, by letting $j\to+\infty$ with the help of the Dominated Convergence Theorem (see \eqref{eq:DeltavjtoDeltaugamma}-\eqref{eq:ToPassLimitDCT}), from \eqref{eq:ugammaDefBd} we obtain
\begin{equation} \label{eq:PDESolvedUgamma}
	u_\gamma(\bar{x},t) = u_0(\bar{x})+\frac{1}{\rho(\bar{x})}\int_0^t \Delta u_\gamma(\bar{x},\tau)\,d\tau\quad\forall\,\,0\leq t\leq T.
\end{equation}
This, together with \eqref{eq:ToUseEquationugamma}, shows that $u_\gamma(\bar{x},\cdot)\in C^1([0,T])$, as desired.
\medskip

\textsc{Step III): $u_\gamma$ solves \eqref{problema}.}
Now we have proved the $C^1$-regularity of $u_\gamma(\bar{x},\cdot)$ (for every fixed $\bar{x}\in G$),
we can now easily show that $u_\gamma$ solves problem \eqref{problema}.
\vspace{0.1cm}

Indeed, let $\bar{x}\in G$ be arbitrarily fixed. On the one hand, since $v_j$ solves problem 
\eqref{eq711}, for every 
$j > d(\bar{x},o)$ we have $v_j(\bar{x},0) = u_0(\bar{x})$;
thus, by letting $j\to+\infty$, from
\eqref{eq:ugammaDefBd} we get
$$u_\gamma(\bar{x},0) = u_0(\bar{x}),$$
and this proves that $u_\gamma$ fulfills the initial condition
in \eqref{problema}. 

On the other hand, owing to \eqref{eq:ToUseEquationugamma}-\eqref{eq:PDESolvedUgamma}, we immediately obtain
\begin{equation} 
	\partial_tu_\gamma(\bar{x},t) = \frac{1}{\rho(\bar{x})}\Delta u_\gamma(\bar{x},t)\quad\forall\,\,0\leq t\leq T,
\end{equation}
and this proves that $u_\gamma$ is a solution of $\LL u = 0$ in $G^T$ (actually, in $G^T_0$).
\medskip

\textsc{Step IV): The condition \eqref{eq:limituNonUnique}.}
Finally, we prove that
 the function $u_\gamma$ (which we
know to be a solution of  \eqref{problema}) satisfies
\eqref{eq:limituNonUnique}. To this end, we fix $0<t_0<T$
and we choose $\varepsilon > 0$ in such a way that
$J = (t_0-\e,t_0+\e)\subseteq (0,T).$
For every $j\in\mathbb{N},\,j > \hat{R}$, we then define
$$w = \gamma-v_j - Ch(x)-\kappa(t-t_0)^2$$
(where $v_j\in\mathfrak{R}^T_0(B_j(o))$ is the
unique solution of the Cauchy\,-\,Di\-ri\-chlet problem \eqref{eq711} introdu\-ced in
the above \textsc{Step I},
and $h$ is as in \eqref{eq710}),
and we claim that
\begin{equation} \label{eq:wleqzeroBarriera}
\text{$w\leq 0$ pointwise on $G\times
\overline{J}$},
\end{equation}
provided that the constants $C,\,\kappa > 0$ are properly chosen.

To prove this claim, it suffices to apply the Weak Maximum Principle
in Lemma \ref{lem:WMP} to the fun\-ction $w$ with the choice
$\Omega = B_j(o)\setminus B_{\hat{R}}(o)$ (see also Remark \ref{rem:zeroImmaterial}). 

Indeed, owing to \eqref{eq710}, and since $v_j$ solves \eqref{eq711} and satisfies
\eqref{eq712}, we have
\begin{align*}
 \mathrm{i)}\,\,&w(x,t_0-\e ) = \gamma-v_j(x,t_0-\e )-Ch(x)-\kappa\,\e^2
 \leq \gamma-\kappa\,\e^2\quad\text{for all $x\in \Omega$}; \\
 \mathrm{ii)}\,\,&w(x,t) = -Ch(x)-\kappa(t-t_0)^2\leq 0\quad\text{for all $x\in G\setminus B_j(o),\,t\in
 \overline{J}$};\\
 \mathrm{iii)}\,\,&w(x,t) \leq \gamma-C\,\min\nolimits_{z\in B_{\hat{R}}(o)}h(z)
 \quad\text{for all $x\in B_{\hat{R}}(o),\,t\in
 \overline{J}$};
 \\
 \mathrm{iv)}\,\,&\LL w = C\Delta h(x)-2\kappa\rho(x)(t-t_0)  \leq
 \rho(x)\big[2\kappa\,\e-C\big]\quad\text{for all $(x,t)\in\Omega\times (t_0-\e,t_0+\e]$}
\end{align*}
 In view of these facts, if we choose $C = C_\e,\,\kappa 
 = \kappa_\e> 0$ in such a way that
 $$1)\,\,\gamma-\kappa\,\e^2\leq 0,\qquad 2)\,\,\gamma-C\,\min_{z\in B_{\hat{R}}(o)}h(z)\leq 0,\qquad
 3)\,\,2\kappa\,\e-C\leq 0$$
 (notice that this is certainly possible, since $h > 0$ pointwise on $G$),
 we are entitled to apply the Weak Maximum Principle
 in Lemma \ref{lem:WMP}, thus obtaining \eqref{eq:wleqzeroBarriera}.
 \medskip

 Now we have established  \eqref{eq:wleqzeroBarriera}, we can easily conclude the proof
 of \eqref{eq:limituNonUnique}. Indeed, owing to the cited  \eqref{eq:wleqzeroBarriera}, and letting
 $j\to+\infty$, by \eqref{eq:ugammaDefBd} we derive that
 $$u_\gamma(x,t) = \lim_{j\to+\infty}v_j(x,t)\geq \gamma-Ch(x)-\kappa(t-t_0)^2\quad
 \text{for all $x\in G,\,t\in J$}.$$
 In particular, choosing $t = t_0\in J$, we get
 $$\gamma-u_\gamma(x,t_0)\leq Ch(x)\quad\text{for all $x\in G$}.$$
From this, since we have already recognized that $\gamma-u_\gamma\geq 0$ on $G^T_0$
 (see again \eqref{eq:ugammaDefBd}), by taking the limit as $d(x,o)\to+\infty$ with the help
 of \eqref{eq710}, we conclude that
 $$\text{$u_\gamma(x,t_0)\to \gamma$ as $d(x,o)\to+\infty$}.$$
 This ends the proof.
\end{proof}

\begin{proof}[Proof of Corollary \ref{cor1pa}] Under the present hypotheses,
it is shown in \cite[Lemma 7.2]{BiMePu}  that there exists a function $h$ as required in Theorem \ref{prop71}; as
a consequence, the thesis follows directly
from Theorem \ref{prop71}.
\end{proof}

\section{Further results on $\mathbb Z^n$: proofs}\label{integer}\setcounter{equation}{0}

We first list two properties of the euclidean distance on the lattice, see also \cite[Theorem 6.1]{MPS1}.
\begin{remark}\label{properties-distance}
Let $x\in\mathbb Z^n$ and consider some $y\in\mathbb Z^n$, $y\sim x$. Then we have, for some $k\in\{1,...,n\}$,
$$
x=(x_1,...,x_n)\quad\text{and}\quad y=(x_1,...,x_k\pm1,...,x_n).
$$
Therefore,
$$
|y|^2-|x|^2=(|x|^2\pm2x_k+1)-|x|^2=\pm2x_k+1\quad\text{and}\quad (|y|^2-|x|^2)^2=4x_k^2+1\pm2x_k\,.
$$
Thus, by summing over all the $y\sim x$ we get
\begin{equation}\label{ae5}
\sum_{y\sim x}\left(|y|^2-|x|^2\right)=2n,\quad\text{and}\quad \sum_{y\sim x}(|y|^2-|x|^2)^2=8|x|^2+2n\,.
\end{equation}
\end{remark}

\begin{proof}[Proof of Theorem \ref{teo3}]
Let us first treat the case \(\alpha\in[0,2)\). Set
\[
\gamma:=\min\{1,2-\alpha\},
\qquad
\beta:=\frac{\gamma}{2}
=\min\left\{\frac12,1-\frac{\alpha}{2}\right\}.
\]
Let \(A>B\) be fixed, and let $Q > 0$.
For \(K:[0, 1/Q]\to(0,+\infty)\), to be specified below, define
\[
\varphi^t(s)=e^{K(t)(1+s)^\beta}
\quad \text{for all } s\geq0.
\]
For any $s, r\in (0, +\infty)$ and for some $\eta$ between $s$ and $r$, we can write
\begin{equation}\label{eq50}
\varphi^t(s)=\varphi^t(r)+(\varphi^t)'(r)(s-r) +\frac{(\varphi^t)''(\eta)}{2}(s-r)^2.
\end{equation}
We compute the derivatives involved in \eqref{eq50},
\begin{align*} 
& (\varphi^t)'(s)=\beta K (1+s)^{\beta -1}\varphi^t(s),\quad(\varphi^t)''(s)\\
& \qquad\quad = \beta(\beta-1)K (1+s)^{\beta-2}\varphi^t(s)+\beta^2K^2(1+s)^{2\beta-2}\varphi^t(s).	
\end{align*}
We now set
$$
K(t):=A(1+Qt),\quad\text{for all}\,\,\,t\in\left[0,\frac1Q\right].
$$
Then we define
$$Z(x,t):=\varphi^t(|x|^2),\quad \text{for all}\,\,\,(x,t)\in\bar S_{\frac1Q}.$$ We now show that, for suitable $A>0$, $Q>0$, $\beta>0$, $Z$ satisfies \eqref{e1p} in 
	$S_{1/Q}$. 

To begin with, we estimate the Laplacian of $Z$. By virtue of \eqref{eq50} with $s=|y|^2,\,r=|x|^2$, we get for all $(x,t)\in \mathbb Z^n\times[0,1/Q]$ with $|x|\geq 2$,
\begin{equation}\label{eq51}
\begin{aligned}
\Delta &Z(x,t)= \frac 1{\mu(x)}\sum_{y\in \mathbb Z^n}[Z(y,t)- Z(x,t)]\omega(x,y) \\
&= \frac 1{2n} \sum_{y\in \mathbb Z^n}\left\{(\varphi^t)'(|x|^2)\left(|y|^2- |x|^2\right)+\frac{(\varphi^t)''(\eta)}{2}\left(|y|^2-|x|^2\right)^2\right\}\omega(x,y)\\
&= \frac 1{2n} \sum_{y\in \mathbb Z^n}\left\{\varphi^t(|x|^2)K(t)\beta(1+|x|^{2})^{\beta-1}(|y|^2- |x|^2)\right.\\
&\left.\quad+\frac{\varphi^t(\eta)}{2}K\beta(1+\eta)^{\beta-2}\left(\beta-1+K(t)\beta(1+\eta)^\beta\right)\left(|y|^2-|x|^2\right)^2\right\}\omega(x,y)\\
&\le \frac {K(t) \beta}{2n} (1+|x|^{2})^{\beta-1}\varphi^t(|x|^2)\sum_{y\sim x}\left(|y|^2- |x|^2\right)\\
&\quad +\frac{K^2(t) \beta^2}{2n}\sum_{y\sim x}\frac{\varphi^t(\eta)}{2}(1+\eta)^{2\beta-2}\left(|y|^2-|x|^2\right)^2\omega(x,y)\\
\end{aligned}
\end{equation}
for some $\eta$ fulfilling
\begin{equation}\label{eq54}
\min\{|x|^2, |y|^2\}\leq \eta\leq \max\{|x|^2, |y|^2\}\,.
\end{equation}
By using \eqref{eq54} and applying the properties of the euclidean distance on the lattice observed in Remark \ref{properties-distance}, we can infer that, if $\beta\in \left(0, \frac 12 \right]$, then
\[\varphi^t(\eta)\leq C \varphi^t(|x|^2)\quad \text{ for all } x\in \mathbb Z^n, |x|\geq 2,\]
for some $C>0$. Therefore, \eqref{eq51} can be furthermore estimated as  \normalcolor
\begin{equation*}
\begin{aligned}
\Delta Z(x,t)&\le K(t) \beta(1+ |x|^{2})^{\beta-1}\varphi^t(|x|^2)\\
& \qquad +C\frac{K^2(t) \beta^2}{2n}\frac{\varphi^t(|x|^2)}{2}(1+|x|^2)^{2\beta-2}\sum_{y\sim x}\left(|y|^2-|x|^2\right)^2\omega(x,y)\\
&\le 2A \beta (1+ |x|^{2})^{\beta-1}\varphi^t(|x|^2)\\
& \qquad +C\frac{4A^2 \beta^2}{2n}\frac{\varphi^t(|x|^2)}{2}(1+|x|^2)^{2\beta-2}\left(8|x|^2+2n\right)\\
&\le 2A \beta (1+ |x|^{2})^{\beta-2}\varphi^t(|x|^2)\times \\
& \qquad \times \left\{(1+ |x|^{2})+C\frac {4A\beta}{n}(1+|x|^{2})^{\beta+1}+AC\beta(1+|x|^{2})^{\beta}\right\}\\
&\le \bar C A^2\beta^2 (1+|x|^2)^{2\beta-1}\varphi^t(|x|^2)\,,
\end{aligned}
\end{equation*}
for a suitable constant 
$$\bar C=\bar C(\beta,C,n,A)>6\max\left\{\frac1{A\beta},\frac{4C\beta}{n},C\right\}.$$ 
Therefore, by means of \eqref{eq:rho-Z} with $\alpha\in[0,2)$, we have
\begin{equation}\label{eq55}
\begin{aligned}
\rho\,\partial_tZ(x,t)&-\Delta Z(x,t)\ge\rho K'(t)Z(x,t)-\bar C A^2 \beta^2 (1+|x|^2)^{2\beta-1}\varphi^t(|x|^2)\\
&\ge \rho_0(1+|x|)^{-\alpha}AQ(1+|x|^2)^{\beta}\varphi^t(|x|^2)-\bar C A^2 \beta^2 (1+|x|^2)^{2\beta-1}\varphi^t(|x|^2)\\
&=(1+|x|^2)^{\beta}\varphi^t(|x|^2)\left\{\rho_0 AQ(1+|x|)^{-\alpha}-\bar C A^2 \beta^2(1+|x|^2)^{\beta-1}\right\}\\
&\ge0\,,
\end{aligned}
\end{equation}
for all $(x,t)\in \mathbb Z^n\times[0,1/Q]$ with
 $|x|\geq 2$, provided that
\[
Q\ge
\frac{\overline C A\beta^2}{\rho_0}
\sup_{|x|\ge2}
(1+|x|)^\alpha(1+|x|^2)^{\beta-1}.
\]
The above supremum is finite, since
\[
\alpha+2\beta-2\le0
\]
by the choice
\[
\beta=\min\left\{\frac12,1-\frac\alpha2\right\}.
\]
Hence, choosing \(Q\) sufficiently large, we get \eqref{eq55}.

\smallskip

On the other hand, the set $\{x\in\mathbb Z^n:\ |x|<2\}$ is finite; moreover, \(\rho(x)>0\) and
\[
\partial_t Z(x,t)
=
AQ(1+|x|^2)^\beta Z(x,t)>0
\]
on this set. Hence, an argument similar to that exploited in
the proof of Theorem \ref{teo1} to obtain
\eqref{eq:DaRichiamareZn}
 shows that, for every $x\in G$ with $|x|<2$ and $0\leq t\leq 1/Q$, 
 \begin{equation}\label{eq58}
\rho(x)\partial_tZ(x,t)-\Delta Z(x,t)\ge0
\end{equation}
by enlarging $Q$ if needed.
Gathering \eqref{eq55} and \eqref{eq58}, we get that $Z$ satisfies \eqref{e1p} in $ S_{1/Q}$. 
\vspace{0.1cm}

We now prove that $Z$ also fulfills \eqref{e3p}, with
$T$ replaced by $T^* = \min\{T,1/Q\}$ (where $Q > 0$
is chosen as above). To this end, according to 
\eqref{eq56}, 
let
\[
\overline{Z}(x):=\exp\left(B|x|^\gamma\right),
\qquad
\gamma=\min\{1,2-\alpha\}.
\]
By assumption,
\[
\limsup_{|x|\to+\infty}
\frac{\max_{t\in[0,T]}u(x,t)}{\overline{Z}(x)}
\le0.
\]
Since $T^*\leq T$, we have
\[
\max_{t\in[0,T^*]}u(x,t)
\le
\max_{t\in[0,T]}u(x,t).
\]
Moreover, by the choice of \(\beta\),
\[
2\beta=\gamma,
\]
and hence
\[
(1+|x|^2)^\beta\sim |x|^\gamma
\qquad\text{as } |x|\to+\infty.
\]
Since \(A>B\), we have, for \(t\in[0,T^*]\),
\[
Z(x,t)
\ge
\exp\left(A(1+|x|^2)^\beta\right),
\]
and therefore
\[
\frac{\overline{Z}(x)}{Z(x,t)}
\le
\frac{\exp\left(B|x|^\gamma\right)}
{\exp\left(A(1+|x|^2)^\beta\right)}
\to0
\qquad\text{as } |x|\to+\infty.
\]
Thus
\[
\limsup_{|x|\to+\infty}
\sup_{t\in[0,T^*]}
\frac{u(x,t)}{Z(x,t)}
\le0.
\]
Hence \(Z\) satisfies condition \eqref{e3p}.
\vspace{0.1cm}

Gathering all these facts, we can exploit
Proposition \ref{prop1}, obtaining
$
u\le0
$
in $\mathbb{Z}^n\times[0,T^*]$.
By ite\-rating this argument, after finitely many steps 
(since \(T<+\infty\) and since \(Q\), being independent of \(T\), can be kept fixed
throughout the iteration) we finally get
\[u\le0
\qquad\text{in }\mathbb{Z}^n\times[0,T]
\]
This is the desired conclusion.
\smallskip

We are left to consider the case \(\alpha=2\). Let again \(A>B\) be fixed, and let  \(Q>0\).
Set
\[
K(t):=A(1+Qt),
\qquad t\in[0,1/Q],
\]
and define
\[
\phi^t(s)=e^{K(t)\log^2(2+s)}
\quad \text{for all }s\ge0.
\]
We compute
$$\begin{aligned}
&(\phi^t)'(s)=\frac{2\,\log(2+s)}{2+s}{K(t)}\phi^t(s),\\
&(\phi^t)''(s)=\frac{4\,\log^2(2+s)}{(2+s)^2}{K^2(t)}\phi^t(s)+2K(t)\left\{\frac{1-\log(2+s)}{(2+s)^2}\right\}\phi^t(s)
\end{aligned}$$
Then, we define
$$Z(x,t):=e^{A(1+Qt)\log^2(2+|x|^2)}=\phi^t(|x|^2),\quad \text{for all}\,\,\,(x,t)\in\bar S_{\frac1Q}.$$
We now show that $Z$ satisfies \eqref{e1p} in 
$S_{1/Q}$. 

To begin with, we estimate the Laplacian of $Z$. By virtue of \eqref{eq50} with $\varphi$ replaced by $\phi$, by choosing $s=|y|^2,\,r=|x|^2$, we get for all $(x,t)\in \mathbb Z^n\times[0,1/Q]$ with $|x|\geq 2$
\begin{equation}\label{eq57}
\begin{aligned}
\Delta &Z(x,t)= \frac 1{\mu(x)}\sum_{y\in \mathbb Z^n}[Z(y,t)- Z(x,t)]\omega(x,y) \\
&= \frac 1{2n} \sum_{y\in \mathbb Z^n}\left\{(\phi^t)'(|x|^2)(|y|^2- |x|^2)+\frac{(\phi^t)''(\eta)}{2}\left(|y|^2-|x|^2\right)^2\right\}\omega(x,y)\\
&= \frac 1{2n} \sum_{y\in \mathbb Z^n}\left\{\phi^t(|x|^2)\frac{2\log(2+|x|^2)}{2+|x|^2}K(t)(|y|^2- |x|^2)\right.\\
&\left.\quad+\frac{\phi^t(\eta)}{2}\frac{K(t)}{(2+\eta)^2}\left[4\log^2(2+\eta)K(t)+2\left(1-\log(2+\eta)\right)\right](|y|^2-|x|^2)^2\right\}\omega(x,y)\\
&\le \frac{K(t)}{n}\frac {\log(2+|x|^2)}{2+|x|^2}\phi^t(|x|^2)\sum_{y\sim x}(|y|^2- |x|^2)\omega(x,y)\\
&\quad +\frac{K(t)^2}{n}\sum_{y\sim x}\frac{\phi^t(\eta)}{(2+\eta)^2}\log^2(2+\eta)\left(|y|^2-|x|^2\right)^2\omega(x,y)\\
\end{aligned}
\end{equation}
for some $\eta$ fulfilling \eqref{eq54}. By using \eqref{eq54}, jointly with the properties of the euclidean distance on the lattice observed in Remark \ref{properties-distance}, \eqref{eq57} can be further estimated, for some $C>0$, with
\begin{equation*}
\begin{aligned}
\Delta Z(x,t)&\le 2K(t)\frac {\log(2+|x|^2)}{2+|x|^2}\phi^t(|x|^2)\\
& \qquad +C\frac{K(t)^2}{n}\phi^t(|x|^2)\frac{\log^2(2+|x|^2)}{(1+|x|^2)^2}\sum_{y\sim x}\left(|y|^2-|x|^2\right)^2\omega(x,y)\\
&\le 4A\frac {\log(2+|x|^2)}{2+|x|^2}\phi^t(|x|^2)+\frac{4A^2}{n}\phi^t(|x|^2)\frac{\log^2(2+|x|^2)}{(2+|x|^2)^2}\left(8|x|^2+2n\right)\\
&\le 4A^2\frac {\log^2(2+|x|^2)}{2+|x|^2}\phi^t(|x|^2)\left\{\frac1{A\log(2+|x|^2)}+\frac {8C}{n}+\frac{2C}{2+|x|^{2}}\right\}\\
&\le \bar C \,A^2\frac {\log^2(2+|x|^2)}{2+|x|^2}\phi^t(|x|^2)\,,
\end{aligned}
\end{equation*}
for some constant 
$$\bar C=\bar C(C,n,A)>12\max\left\{\frac1{A\log3},\frac{8C}{n},\frac{2C}{3}\right\}.$$ 
Therefore, by means of \eqref{eq:rho-Z} with $\alpha=2$, we have 
\begin{equation}\label{eq510}
\begin{aligned}
\rho\,\partial_tZ(x,t)&-\Delta Z(x,t)\ge\rho K'(t)Z(x,t)-\bar C \,A^2\frac {\log^2(2+|x|^2)}{2+|x|^2}\phi^t(|x|^2)\\
&\ge \frac{\rho_0}{4(1+|x|^{2})}AQ\log^2(2+|x|^2)\phi^t(|x|^2)-\bar C \,A^2\frac {\log^2(2+|x|^2)}{2+|x|^2}\phi^t(|x|^2)\\
&=A\frac{\log^2(2+|x|^2)}{1+|x|^{2}}\phi^t(|x|^2)\left\{\frac{\rho_0}{4} Q-\bar C A\right\}\\
&\ge0\,,
\end{aligned}
\end{equation}
for all $(x,t)\in \mathbb Z^n\times[0,1/Q]$
with $|x|\geq 2$, provided that $Q\ge \frac{4\bar CA}{\rho_0}$.
\vspace{0.1cm}

On the other hand, the set $\{x\in\mathbb Z^n:\ |x|<2\}$ is finite. Moreover, \(\rho(x)>0\) and
\[
\partial_tZ(x,t)
=
AQ\log^2(2+|x|^2)Z(x,t)>0
\]
on this set. Hence, by arguing as above
(that is, with an argument similar to that exploited
in the proof of Theorem 
\ref{teo1}), 
for every \(|x|<2\) and every \(t\in[0,1/Q]\) we get
\begin{equation}\label{eq511}
\rho(x)\partial_tZ(x,t)-\Delta Z(x,t)\ge0
\end{equation}
by enlarging $Q$ if needed.
Gathering \eqref{eq510} and \eqref{eq511}, we get that $Z$ satisfies \eqref{e1p} in $S_{1/Q}$. 
\vspace{0.1cm}

We now prove that $Z$ also fulfills \eqref{e3p}, with
$T$ replaced by $T^* = \min\{T,1/Q\}$ (where $Q > 0$
is chosen as above). To this end, according
to \eqref{eq56}, set
\[
\overline{Z}(x):=\exp\left(B\log^2(2+|x|^2)\right).
\]
By assumption,
\[
\limsup_{|x|\to+\infty}
\frac{\max_{t\in[0,T]}u(x,t)}{\overline{Z}(x)}
\le0.
\]
Since \(T^*\le T\), we have
\[
\max_{t\in[0,T^*]}u(x,t)
\le
\max_{t\in[0,T]}u(x,t).
\]
Moreover, since \(A>B\), for every \(t\in[0,T^*]\),
\[
Z(x,t)
\ge
\exp\left(A\log^2(2+|x|^2)\right).
\]
Therefore
\[
\frac{\overline{Z}(x)}{Z(x,t)}
\le
\exp\left(-(A-B)\log^2(2+|x|^2)\right)
\to0
\qquad\text{as } |x|\to+\infty.
\]
Thus
\[
\limsup_{|x|\to+\infty}
\sup_{t\in[0,1/Q]}
\frac{u(x,t)}{Z(x,t)}
\le0.
\]
Hence \(Z\) satisfies condition \eqref{e3p}.
\vspace{0.1cm}

Gathering all these facts, we can exploit
Proposition \ref{prop1}, obtaining
$
u\le0
$
in $\mathbb{Z}^n\times[0,T^*]$.
By ite\-rating this argument, after finitely many steps 
(since \(T<+\infty\) and since \(Q\), being independent of \(T\), can be kept fixed
throughout the iteration) we finally get
\[u\le0
\qquad\text{in }\mathbb{Z}^n\times[0,T]
\]
This ends the proof.
\end{proof}

\normalcolor 
\begin{proof}[Proof of Corollary \ref{cor3pa}] Under the present hypotheses in the proof of \cite[Proposition 4.3]{BP25} it is shown that there exists a function $h$ as required in Theorem \ref{prop71}.
Hence the thesis follows from Theorem \ref{prop71}.
\end{proof}

\section{Further results on antitrees}\label{antitree}\setcounter{equation}{0}

In this section, we state uniqueness and non-uniqueness results on anti-trees. In particular, we see that in general the uniqueness condition of Theorem \ref{teo1} is not sharp in this context. However, the new condition that we find on anti-trees is optimal. 

\subsection{Anti-trees and radial functions} We keep the notation as in Subsection \ref{degree}. Let $\Omega= \{o\}$ for some point $o\in G$.
Let $s: \mathbb N\to \mathbb N$ be given by
$$s(r) = \operatorname{card}[S_r(o)] \quad \textrm{for all } r \in \mathbb N. $$
We say that \(G\) is an anti-tree with sphere size \(s\) if the
following properties hold.  Each vertex in \(S_r(o)\), \(r\geq 1\),
is adjacent to every vertex in \(S_{r-1}(o)\) and to every vertex in
\(S_{r+1}(o)\).  Moreover, there are no edges between vertices
belonging to the same sphere.
By the definition of anti-tree, one has
\[
   \mathfrak  D_+(x)=s(r+1),
    \qquad
    \mathfrak D_-(x)=s(r-1),
    \qquad x\in S_r(o),\ r\geq 1.
\]
At the root \(o\), we have
\[
    \mathfrak D_+(o)=s(1),
    \qquad
    \mathfrak D_-(o)=0 .
\]
Let \(F:\mathbb N_0\to\mathbb R\), and let \(f:V\to\mathbb R\) be the
radial function defined by
\[
    f(x)=F(r),
    \qquad r=d(x,o).
\]
If \(x\in S_r(o)\), \(r\geq 1\), then the neighbours of \(x\) belong
only to \(S_{r-1}(o)\) and \(S_{r+1}(o)\).  In view of \eqref{e13p}, we have 
\begin{equation}\label{enf1}
  \Delta f(x)
    =
    s(r+1)\bigl(F(r+1)-F(r)\bigr)
    +
    s(r-1)\bigl(F(r-1)-F(r)\bigr)
\quad \text{ for any } x\in S_r(o), r\geq 1.  
\end{equation}
At the root one has instead
\begin{equation}\label{enf1b}
    \Delta f(o)=s(1)\bigl(F(1)-F(0)\bigr).
\end{equation}


 \subsection{Uniqueness on antitrees}
\begin{lemma}\label{alemma1} 
Let \(\rho:G\to(0,+\infty)\). Assume that there exists a positive
sequence \(\{\rho_r\}_{r\in \mathbb N_0}\) such that
\[
    \rho(x)\geq \rho_r
    \qquad \text{for every } x\in S_r(o).
\]
Assume moreover that
\begin{equation}\label{enf2}
    \sum_{n=1}^{+\infty}
     \frac{\sum_{k=0}^{n-1}s(k)\rho_k}
         {s(n-1)s(n)}
    =
    +\infty .
\end{equation}
Then there exists a radial function \(Z:G\to(0,+\infty)\) such that
\begin{equation}\label{enf3}
    Z(r)\to+\infty
    \qquad \text{as } r\to+\infty,
\end{equation}
and
\begin{equation}\label{enf4}
    \Delta Z(x)\leq \rho(x)
    \qquad \text{for every } x\in G.
\end{equation}
\end{lemma}

\begin{proof}
Let \(M>0\). We define \(Z:G\to(0,+\infty)\) as a radial function. Namely,
we set
\[
    Z(o)=Z(0):=M,
\]
and, for \(x\in S_r(o)\), \(r\geq1\),
\begin{equation}\label{enf7}
    Z(x)=Z(r):=
    M+
    \sum_{n=1}^{r}
    \frac{\sum_{k=0}^{n-1}s(k)\rho_k}
         {s(n-1)s(n)} .
\end{equation}
Equivalently, \eqref{enf7} also holds for \(r=0\), with the convention that
the empty sum is equal to zero.

By construction, \(Z\geq M>0\). Moreover, \eqref{enf2} gives \eqref{enf3}.
It remains to check \eqref{enf4}. Set
\[
    b(r):=Z(r)-Z(r-1),
    \qquad r\geq1.
\]
By the definition of \(Z\),
\[
    b(r)=
    \frac{\sum_{k=0}^{r-1}s(k)\rho_k}
         {s(r-1)s(r)} .
\]
Let
\[
    A(r):=s(r-1)s(r)b(r),
    \qquad r\geq1.
\]
Then
\[
    A(r)=\sum_{k=0}^{r-1}s(k)\rho_k .
\]
In particular,
\[
    A(1)=s(0)\rho_0=\rho_0,
\]
because \(s(0)=1\), and
\[
    A(r+1)-A(r)=s(r)\rho_r,
    \qquad r\geq1.
\]

In view of \eqref{enf1b}, we get
\[
    \Delta Z(o)=s(1)b(1)=A(1)=\rho_0.
\]
Hence
\begin{equation}\label{enf5}
    \Delta Z(o)=\rho_0\leq \rho(o).
\end{equation}

Let now \(x\in S_r(o)\), with \(r\geq1\). By \eqref{enf1},
\[
    \Delta Z(x)=\frac{A(r+1)-A(r)}{s(r)}.
\]
Therefore
\begin{equation}\label{enf6}
    \Delta Z(x)=\frac{s(r)\rho_r}{s(r)}=\rho_r\leq \rho(x)
    \qquad \text{for every } x\in S_r(o),\ r\geq1.
\end{equation}
From \eqref{enf5} and \eqref{enf6} we deduce \eqref{enf4}. This completes
the proof.
\end{proof}

From Lemma \ref{lemmaelliptic} and Proposition \ref{prop1}, we can immediately deduce the following
\begin{theorem}\label{teoa1}
Let the assumptions of Lemma \ref{alemma1} be satisfied. Let $u$ be a subsolution of problem \eqref{problema} with $f\equiv u_0\equiv0$ fulfilling
\begin{equation*}
\lim_{r\to +\infty}\frac1{Z(r)}\left\{\max_{t\in[0,T]}{|u(x,t)|}\right\} = 0\,,
\end{equation*}
with $Z$ defined in \eqref{enf7}. Then
\[u \leq 0 \quad\text{in}\,\,\, S_T.\]
\end{theorem}

\begin{corollary}\label{cora3Bis}
Let the assumptions of Lemma \ref{alemma1} be satisfied.  Then there exists at most one solution to problem \eqref{problema} such that
\[\lim_{r\to +\infty}\frac1{Z(r)}\left\{\max_{t\in[0,T]}{|u(x,t)|}\right\} =0\,,\]
with $Z$ defined in \eqref{enf7}. 
\end{corollary}

\begin{remark}
In the special case \(\rho\equiv1\), and when only bounded solutions are
considered, Corollary~\ref{cora3Bis} is in agreement with
\cite[Example~9.28]{KLW}, in which, under the same assumptions, the graph is
shown there to be stochastically complete. Hence problem \eqref{problema}
admits at most one bounded solution.
\end{remark}

\begin{example}
We now give a simple example showing that the uniqueness criterion in Theorem \ref{teoa1} 
may apply even when the general Theorem \ref{teo1} cannot be used.  Let \(G\) be the anti-tree with sphere size
\[
    s(r)=r+1,
    \qquad r\geq 0 .
\]
Then, for every \(x\in S_r(o)\), \(r\geq 1\), one has
\[
    \mathfrak D_+(x)=s(r+1)=r+2.
\]
Hence
\[
    \frac{\mathfrak D_+(x)}{r}
    =
    \frac{r+2}{r}
    \sim 1
    \qquad \text{as } r\to+\infty .
\]
Let now
\[
    \rho(x)=\rho_r:=(r+2)^{-1/2},
    \qquad x\in S_r(o).
\]
Then
\[
    \rho_r=o\left(\frac{\mathfrak D_+(x)}{r}\right)
    \qquad \text{as } r\to+\infty,
\]
because
\[
    \frac{\rho_r}{\mathfrak D_+(x)/r}
    =
    \frac{r}{(r+2)^{3/2}}
    \to 0 .
\]
Therefore, Theorem \ref{teo1} cannot be applied. Nevertheless, \eqref{enf2} holds.
Indeed,
\[
\begin{aligned}
    \sum_{n=1}^{+\infty}
    \frac{\sum_{k=0}^{n-1}s(k)\rho_k}
         {s(n-1)s(n)}
    &=
    \sum_{n=1}^{+\infty}
    \frac{\sum_{k=0}^{n-1}(k+1)(k+2)^{-1/2}}
         {n(n+1)} .
\end{aligned}
\]
For \(k\in\{\lfloor n/2\rfloor,\ldots,n-1\}\), we have
\[
    k+1\geq \frac{n}{2},
    \qquad
    (k+2)^{-1/2}\geq C n^{-1/2},
\]
for a positive constant \(C\) independent of \(n\).  Therefore
\[
    \sum_{k=0}^{n-1}(k+1)(k+2)^{-1/2}
    \geq C n^{3/2}.
\]
Consequently,
\[
    \frac{\sum_{k=0}^{n-1}(k+1)(k+2)^{-1/2}}
         {n(n+1)}
    \geq C n^{-1/2}.
\]
Since
\[
    \sum_{n=1}^{+\infty} n^{-1/2}=+\infty,
\]
we obtain \eqref{enf2}. Thus, Theorem \ref{teoa1} can be applied. 
\end{example}

\subsection{Non-uniqueness on antitrees}
\begin{lemma}\label{h-antit}
Let \(G\) be an anti-tree with size \(s\). Let \(\rho:G\to(0,+\infty)\). 
Assume that there exists a positive sequence \(\{\rho_r\}_{r\in \mathbb N_0}\) 
such that
\[
    \rho(x)\leq \rho_r
    \qquad \text{for every } x\in S_r(o).
\]
Assume moreover that
\begin{equation}\label{enf8}
    \sum_{n=1}^{+\infty}
    \frac{\sum_{k=0}^{n-1}s(k)\rho_k}
         {s(n-1)s(n)}
    <
    +\infty .
\end{equation}
Then there exists a radial function \(h:G\to(0,+\infty)\) such that
\begin{equation}\label{enf9}
    h(r)\to0
    \qquad \text{as } r\to+\infty,
\end{equation}
and
\begin{equation}\label{enf10}
    \Delta h(x)\leq -\rho(x)
    \qquad \text{for every } x\in G.
\end{equation}
\end{lemma}

\begin{proof}
We define \(h:G\to(0,+\infty)\) as a radial function. Namely, for
\(x\in S_r(o)\), we set
\[
    h(x)=h(r),
\]
where
\[
    h(r):=
    \sum_{n=r+1}^{+\infty}
    \frac{\sum_{k=0}^{n-1}s(k)\rho_k}
         {s(n-1)s(n)} ,
    \qquad r\geq 0.
\]
The series is finite by assumption \eqref{enf8}. Since all its terms are
positive, we have
\[
    h(r)>0
    \qquad \text{for every } r\geq 0.
\]
Moreover, \(h(r)\) is the tail of a convergent positive series, hence
\eqref{enf9} holds.

It remains to check \eqref{enf10}. Set
\[
    a(r):=h(r-1)-h(r),
    \qquad r\geq 1.
\]
By the definition of \(h\),
\[
    a(r)=
    \frac{\sum_{k=0}^{r-1}s(k)\rho_k}
         {s(r-1)s(r)} .
\]
Let
\[
    A(r):=s(r-1)s(r)a(r),
    \qquad r\geq 1.
\]
Then
\[
    A(r)=\sum_{k=0}^{r-1}s(k)\rho_k .
\]
In particular,
\[
    A(1)=s(0)\rho_0=\rho_0,
\]
because \(s(0)=1\), and
\[
    A(r+1)-A(r)=s(r)\rho_r,
    \qquad r\geq 1.
\]

From \eqref{enf1b} we get
\[
    \Delta h(o)
    =
    s(1)\bigl(h(1)-h(0)\bigr)
    =
    -s(1)a(1).
\]
Since
\[
    A(1)=s(0)s(1)a(1)=s(1)a(1),
\]
we obtain
\[
    \Delta h(o)=-A(1)=-\rho_0.
\]
By the assumption \(\rho(o)\leq \rho_0\), it follows that
\begin{equation}\label{enf12}
    \Delta h(o)\leq -\rho(o).
\end{equation}

Let now \(x\in S_r(o)\), with \(r\geq 1\). Since
\[
    h(r+1)-h(r)=-a(r+1),
    \qquad
    h(r-1)-h(r)=a(r),
\]
\eqref{enf1} gives
\[
\begin{aligned}
    \Delta h(x)
    &=
    -s(r+1)a(r+1)+s(r-1)a(r)       \\
    &=
    -\frac{A(r+1)}{s(r)}
    +
    \frac{A(r)}{s(r)}                \\
    &=
    -\frac{A(r+1)-A(r)}{s(r)}.
\end{aligned}
\]
Using \(A(r+1)-A(r)=s(r)\rho_r\), we get
\begin{equation}\label{enf13}
    \Delta h(x)=-\rho_r\leq -\rho(x)
    \qquad \text{for every } x\in S_r(o),\ r\geq 1.
\end{equation}
From \eqref{enf12} and \eqref{enf13} we obtain \eqref{enf10}. This completes
the proof.
\end{proof}

From Lemma \ref{h-antit} and Theorem \ref{prop71} we obtain the next non-uniqueness result.
\begin{theorem}\label{teo2antit}
Let the assumptions of Lemma \ref{h-antit}  be satisfied.  Then there exist infinitely many bounded solutions $u$ of problem \eqref{problema}. In
particular, for every fixed $\gamma\in \R$ fulfilling 
\eqref{eq:assumptionOnuzero} 
there exists a solution $u$ to problem \eqref{problema} such that
\eqref{eq:limituNonUnique} holds.
\end{theorem}

\section{The special case of $\mathbb Z^2$}\label{Z2}\setcounter{equation}{0}

In this section, we demonstrate that on $\mathbb Z^2$ problem \eqref{problema} admits a unique solution satisfying an appropriate growth condition at infinity, for every
$\rho\in \mathfrak{F} , \rho>0$ in $G$.
This reveals a striking contrast between the behavior of $\mathbb Z^2$ and the cases previously examined in Sections \ref{main} and \ref{Zn}.

\begin{lemma}\label{alemma3}
Let $K>0$ and
\begin{equation}\label{ae23}
\hat Z(x):= K \log(\log(|x|^2+4)) \quad \text{ for any } x\in \mathbb Z^2.
\end{equation}
Then, for some $K>0$,
\begin{equation}\label{ae20}
\Delta \hat Z \leq \rho(x) \quad \text{ for any } x\in \mathbb Z^2\,.
\end{equation}
\end{lemma}
\begin{proof}
 The proof of this lemma is entirely based on following key fact concerning the function $\hat Z$
 defined in \eqref{ae23}: \emph{it is possible to find a positive
 number $R_0 > 0$ such that}
 \begin{equation} \label{eq:claimDeltalog}
  \Delta\big(z\mapsto \log(\log(4+|z|^2))\big)(x)< 0\quad\text{for every $x\in\mathbb{Z}^2,\,|x| > R_0$}.
 \end{equation}
 Taking this fact for granted for a moment, we can easily prove \eqref{ae20}. Indeed,
 let $R_0 > 0$ be as in \eqref{eq:claimDeltalog}. Since the ball $B_{R_0}(0)$ is a finite set,
 and since $\rho > 0$ pointwise on $G$, we have
 \begin{equation} \label{eq:deltaHatZVicino}
  \begin{split}
  \Delta \hat{Z}(x) & = K\cdot\Delta\big(z\mapsto \log(\log(4+|z|^2))\big)(x)\\
  & \leq K\cdot\max_{B_{R_0}(0)}\big|\Delta\big(z\mapsto \log(\log(4+|z|^2))\big)(\cdot)\big|
  \\
  & \leq \min_{B_{R_0}(0)}\rho(\cdot)\leq \rho(x)\quad\text{for every $x\in B_{R_0}(0)$},
  \end{split}
 \end{equation}
 provided that $K > 0$ is small enough. On the other hand,
 by \eqref{eq:claimDeltalog} we also have
 \begin{equation} \label{eq:deltaHatZLontano}
  \begin{split}
   \Delta \hat{Z}(x) & = K\cdot\Delta\big(z\mapsto \log(\log(4+|z|^2))\big)(x)
   \\
   & <  0 < \rho(x)\quad\text{for every $x\in\mathbb{Z}^2,\,|x| > R_0$}.
  \end{split}
 \end{equation}
 Thus, by combining \eqref{eq:deltaHatZVicino}\,-\,\eqref{eq:deltaHatZLontano} we immediately obtain
 \eqref{ae20}.
 \vspace*{0.1cm}

 Hence, we turn to prove \eqref{eq:claimDeltalog}.  To this end we first observe that,
 setting
 $$\varphi(t) = \log(\log(4+t)),$$
 by using the Taylor formula with Lagrange remained
 (and by taking into account
 the explicit expression of $\omega$ and of $\mu$ in this setting,
 see Section \ref{Zn}), for every $x\in\mathbb{Z}^2$ we can write
 \begin{align*}
  & \Delta\big(z\mapsto \log(\log(4+|z|^2))\big)(x)
  = \frac{1}{4}\sum_{y\sim x}\big[\varphi(|y|^2)-\varphi(|x|^2)\big]
  \\
  & \qquad\qquad= \frac{1}{4}\Big\{
  \varphi'(|x|^2)\sum_{y\sim x}(|y|^2-|x|^2)
  +\frac{\varphi''(|x|^2)}{2}\sum_{y\sim x}(|y|^2-|x|^2)^2 \\
 &\qquad\qquad\qquad\qquad
 +\sum_{y\sim x}\frac{\varphi^{(3)}(\xi_{x,y})}{6}(|y|^2-|x|^2)^3\Big]\Big\} = (\bigstar),
 \end{align*}
 where $\xi_{x,y} \in\R$ is a  point between $|x|^2$ and $|y|^2$;
 from this, by \eqref{ae5} we obtain
 \begin{equation} \label{eq:DeltaLogdaStim}
 \begin{split}
  (\bigstar) & = \frac{1}{4}\Big\{
  4\varphi'(|x|^2)+(4|x|^2+2)\varphi''(|x|^2)
   +\sum_{y\sim x}\frac{\varphi^{(3)}(\xi_{x,y})}{6}(|y|^2-|x|^2)^3\Big\} \\
   & = \frac{1}{4}\big\{A(x)+B(x)\big\},
 \end{split}
 \end{equation}
 where we have introduce the notation
 \begin{align*}
  A(x) & = 4\varphi'(|x|^2)+(4|x|^2+2)\varphi''(|x|^2) \\[0.15cm]
   B(x) & =
  \sum_{y\sim x}\frac{\varphi^{(3)}(\xi_{x,y})}{6}(|y|^2-|x|^2)^3
 \end{align*}
 We then estimate the two terms $A(x),\,B(x)$ as $|x|\to+\infty$.
 \medskip

 \noindent-\,\,\emph{Estimate of $A(x)$}. By explicitly computing the derivatives of $\varphi$, we get
 \begin{align*}
  A(x) & = \frac{4}{(4+|x|^2)\log(4+|x|^2)}
  -\frac{(4|x|^2+2)(1+\log(4+|x|^2))}{(4+|x|^2)^2\log^2(4+|x|^2)} \\
  & = \frac{1}{(4+|x|^2)^2\log^2(4+|x|^2)}\times \\
  & \qquad\qquad\times\big[4(4+|x|^2)\log(4+|x|^2)-(4|x|^2+2)(1+\log(4+|x|^2))\big] \\
  & \sim \frac{1}{2|x|^4\log^2(|x|)}\cdot(-4|x|^2) = -\frac{2}{|x|^2\log^2(|x|)}\qquad
  \text{as $|x|\to+\infty$};
 \end{align*}
 as a consequence, there exists $R_0 > 0$ such that
 \begin{equation} \label{eq:estimAsintoticaA}
  A(x) \leq -\frac{1}{|x|^2\log^2(|x|)}\quad
  \text{for every $x\in\mathbb{Z}^2$ with $|x| > R_0$}.
 \end{equation}
 -\,\,\emph{Estimate of $B(x)$}. First of all we observe that,
 by computing $\varphi^{(3)}(t)$, we have
 \begin{equation*}
 \begin{split}
  \varphi^{(3)}(t)
   & = \frac{2(1+\log^2(t+4))(t+4)\log(t+4)-(t+4)\log^2(t+4)}{(t+4)^4\log^4(t+4)} \\
   & \sim \frac{2(t+4)\log^3(t+4)}{(t+4)^4\log^4(t+4)}
   \sim \frac{2}{t^3\log(t)}\qquad\text{as $t\to+\infty$};
  \end{split}
 \end{equation*}
 thus, we can find some $t_0 > 0$ such that
 \begin{equation} \label{eq:stimaDerTerza}
  0\leq \varphi^{(3)}(t) \leq \frac{4}{t^3\log(t)}\qquad\text{for every $t > t_0$}.
 \end{equation}
 On the other hand, if $y\in\mathbb{Z}^2$ and if $y\sim x$,
 we have $y = x\pm e_i$, where $e_i$ is the $i$-th vector
 of the canonical basis of $\R^2$ (for $i = 1,2$); hence, for every $x\in\mathbb{Z}^2$ with
 $|x|\geq 2$ we get
 $$(|x|-1)^2\leq |y|^2\leq (|x|+1)^2,$$
 From this, since $\xi_{x,y}$ is between $|x|^2$ and $|y|^2$, we deduce that
 \begin{equation}  \label{eq:estimaxixy}
 \begin{gathered}
  \text{$(|x|-1)^2\leq \xi_{x,y}\leq (|x|+1)^2$ } \\
  \text{for every $x\in\mathbb{Z}^2$ with $|x|\geq 2$ and every $y\sim x$}.
  \end{gathered}
  \end{equation}
 Summing up, by combining \eqref{eq:stimaDerTerza}\,-\,\eqref{eq:estimaxixy}
 (and by possibly
 enlarging the number $R_0$ introdu\-ced in \eqref{eq:estimAsintoticaA} in such a way
 that $(R_0-1)^2\geq t_0$), we obtain
 \begin{equation*}
  \begin{split}
   B(x) & \leq
   \frac{2}{(|x|-1)^6\log(|x|-1)}\cdot\frac{1}{6}\sum_{y\sim x}|(|y|^2-|x|^2)^3| \\
   & \leq  \frac{4(1+2|x|)^3}{3(|x|-1)^6\log(|x|-1)}
   \qquad\text{for every $x\in\mathbb{Z}^2$ with $|x| > R_0$}.
  \end{split}
 \end{equation*}
 This, together with the obvious asymptotic equivalence
 $$\frac{4(1+2|x|)^3}{3(|x|-1)^6\log(|x|-1)}
 \sim \frac{32}{3|x|^3\log(|x|)}\quad\text{as $|x|\to+\infty$},
 $$
 finally gives the following estimate for $B(x)$
 \begin{equation} \label{eq:estimAsintoticaB}
  B(x)\leq \frac{32}{|x|^3\log(|x|)}\quad\text{for every $x\in\mathbb{Z}^2$ with $|x| > R_0$}
 \end{equation}
  (up to possibly enlarging once again the number $R_0$).
  \medskip

 Now we have estimated the terms $A(x)$ and $B(x)$, we can easily conclude
 the demonstration
 of the claimed \eqref{eq:claimDeltalog}: indeed,
 by combining \eqref{eq:estimAsintoticaA} with
 \eqref{eq:estimAsintoticaB}, from \eqref{eq:DeltaLogdaStim} we get
 \begin{align*}
  & \Delta\big(z\mapsto \log(\log(4+|z|^2))\big)(x)
  \leq \frac{1}{4}\big\{A(x)+B(x)\big\} \\
  & \qquad\leq -\frac{1}{2|x|^2\log^2(|x|)}
  + \frac{8}{|x|^3\log(|x|)} \\
  & \qquad =
  -\frac{1}{2|x|^2\log^2(|x|)}\Big(1-\frac{16\log(|x|)}{|x|}\Big);
 \end{align*}
 as a consequence, since we clearly have
 $$1-\frac{16\log(|x|)}{|x|}\to 1\quad\text{as $|x|\to+\infty$},$$
 by possibly enlarging $R_0 > 0$ we conclude that
 $$\Delta\big(z\mapsto \log(\log(4+|z|^2))\big)(x)
 \leq -\frac{1}{4|x|^2\log^2(|x|)}< 0,$$
 for every $x\in\mathbb{Z}^2$ with $|x| > R_0$. This ends the proof.
\end{proof}


From Lemma \ref{alemma3} and Proposition \ref{prop1}, we can immediately deduce the following
\begin{theorem}\label{teoa3}
Let $\rho\in \mathfrak F, \rho>0$ in $\mathbb Z^2$. Let $u$ be a subsolution of problem \eqref{problema} with $f\equiv u_0\equiv0$ fulfilling
\begin{equation*}
\lim_{|x|\to +\infty}\frac1{\log(\log|x|^2))}\left\{\max_{t\in[0,T]}{|u(x,t)|}\right\} = 0\,,
\end{equation*}
Then
\[u \leq 0 \quad\text{in}\,\,\, S_T.\]
\end{theorem}

\begin{corollary}\label{cora3}
Let $\rho\in \mathfrak F, \rho>0$ in $\mathbb Z^2$. Then there exists at most one solution to problem \eqref{problema} such that
\[\lim_{|x|\to +\infty} \frac1{\log(\log|x|^2))}\left\{\max_{t\in[0,T]}{|u(x,t)|}\right\} =0\,.\]
\end{corollary}

\bigskip
\bigskip

\noindent{\bf Acknowledgement}.
All authors are member of the ``Gruppo Nazionale per l'Analisi Ma\-te\-ma\-tica, la Probabilit\`a e le loro Applicazioni'' (GNAMPA) of the ``Istituto Nazionale di Alta Matematica'' (INdAM, Italy).
The first author is partially supported by the PRIN 2022 project 2022R537CS \emph{$NO^3$ - Nodal Optimization, NOnlinear elliptic equations, NOnlocal geometric problems, with a focus on regularity}, founded by the European Union - Next Generation EU.
The second author is funded by the Deutsche Forschungsgemeinschaft (DFG, German Research Foundation) - SFB 1283/2 2021 - 317210226.
The third author acknowledge that this work is part of the PRIN project 2022 Geometric-analytic methods for PDEs and applications, ref. 2022SLTHCE, financially supported by the EU, in the framework of the "Next Generation EU initiative".

\end{document}